\documentclass[10pt,journal]{IEEEtran}
\ifCLASSOPTIONcompsoc
  \usepackage[nocompress]{cite}
\else
  \usepackage{cite}
\fi
\usepackage{amssymb}
\usepackage{amsmath}
\usepackage{dsfont}
\usepackage{amsthm}
\usepackage{graphicx}
\usepackage{mathrsfs,doi,color}
\usepackage{subfigure}

\newcommand{\real}{\mathbb{R}}
\newcommand{\complex}{\mathbb{C}}

\newtheoremstyle{mythm}{1.5ex plus 1ex minus .2ex}{1.5ex plus 1ex minus .2ex} {\rm}{\parindent}{\it\it}{\rm{:}}{1em}{}
\theoremstyle{mythm}

\renewcommand{\IEEEQED}{\IEEEQEDopen}

\newtheorem{theorem}{Theorem}[section]
\newtheorem{corollary}{Corollary}[section]
\newtheorem{lemma}{Lemma}[section]
\newtheorem{proposition}{Proposition}[section]
\newtheorem{definition}{Definition}[section]
\newtheorem{remark}{Remark}

\newcommand{\E}{\mathbb{E}}

\ifCLASSINFOpdf
\else
\fi
\begin{document}

\title{Linear Stochastic Approximation Algorithms and Group Consensus over Random Signed Networks: A Technical Report with All Proofs}

\author{Ge~Chen,~\IEEEmembership{Member,~IEEE}, Xiaoming Duan, Wenjun Mei,  Francesco Bullo,~\IEEEmembership{Fellow,~IEEE}
  \IEEEcompsocitemizethanks{\IEEEcompsocthanksitem This material is
    based upon work supported by, or in part by, the U.S.\ Army
    Research Laboratory and the U.S.\ Army Research Office under grant
    number W911NF-15-1-0577. The research of G.\ Chen was supported in part
by the National Natural Science Foundation of China under grants 91427304,
61673373 and 11688101, the National Key Basic Research Program of China
(973 program) under grant 2014CB845301/2/3, and the Leading research
projects of Chinese Academy of Sciences under grant QYZDJ-SSW-JSC003.

  \IEEEcompsocthanksitem
  Ge Chen is with the National Center for Mathematics and Interdisciplinary Sciences \& Key Laboratory of Systems and
  Control, Academy of Mathematics and Systems Science, Chinese Academy of Sciences, Beijing 100190,
  China, {\tt  chenge@amss.ac.cn}
  \IEEEcompsocthanksitem Xiaoming Duan, Wenjun Mei, and Francesco Bullo
  are with the Department of Mechanical Engineering and the Center of
  Control, Dynamical-Systems and Computation, University of California
  at Santa Barbara, CA 93106-5070, USA. {\tt
    xiaomingduan.zju@gmail.com; meiwenjunbd@gmail.com;
    bullo@engineering.ucsb.edu}
}
}

\IEEEtitleabstractindextext{%
  \begin{abstract}
    This paper studies linear stochastic approximation (SA) algorithms
    and their application to multi-agent systems in engineering and
    sociology.
    As main contribution, we provide necessary and sufficient
    conditions for convergence of linear SA algorithms to a deterministic or random final vector.  We
    also characterize the system convergence rate, when the system is
    convergent. Moreover, differing from non-negative gain functions
    in traditional SA algorithms, this paper considers also the case
    when the gain functions are allowed to take arbitrary real
    numbers.
    Using our general treatment, we provide necessary and sufficient
    conditions to reach consensus and group consensus for first-order
    discrete-time multi-agent system over random signed networks and
    with state-dependent noise.
    Finally, we extend our results to the setting of multi-dimensional
    linear SA algorithms and characterize the behavior of the
    multi-dimensional Friedkin-Johnsen model over random interaction
    networks.
\end{abstract}

}

\maketitle

\IEEEdisplaynontitleabstractindextext

\IEEEpeerreviewmaketitle

%
%
\renewcommand{\thesection}{\Roman{section}}
\section{Introduction}
\renewcommand{\thesection}{\arabic{section}}
Distributed coordination of multi-agent systems has drawn much
attention from various fields over the past decades. For example,
engineers control the formations of mobile robots, satellites,
unmanned aircraft, and automated highway systems
\cite{JAF-RMM:04,WR-RWB:08}; physicists and computer scientists model
the collective behavior of animals \cite{TV-AC-EBJ-IC-OS:95,CWR:87};
sociologists investigate the evolution of opinion, belief and social
power over social networks
\cite{MHDG:74,PJ-AM-NEF-FB:13d,NEF-AVP-RT-SEP:16}.  Many models for
distributed coordination have been proposed and analyzed; a common
thread in all these works is the study of a group of interacting
agents trying to achieve a collective behavior by using neighborhood
information allowed by the network topology.

Linear dynamical systems are a class of basic first-order dynamics with
application to many practical problems in multi-agent systems,
including distributed consensus of multi-agent systems, computation of
PageRank, sensor localization of wireless networks, opinion dynamics, and belief
evolution on social networks
\cite{CR-PF-RT-HI:15,AVP-RT:17,NEF-AVP-RT-SEP:16}. If the operator in a linear dynamical system is time-invariant, then the study of this system is
straightforward.
However, practical systems are very often subject to random
fluctuations, so that the operator in an linear dynamical system is
time-variant and the system may not converge.  To overcome
this deficiency and eliminate the effects of fluctuation, a feasible
approach is to adopt models based on the stochastic approximation (SA)
algorithm \cite{MYH-JHM:09,RC-GC-PF-FG:11,NEL-AO:14,TL-JFZ:10,MH:12,HT-TL:15,GL-CG:17}.

The main idea of the SA algorithm is as follows: each agent has a
memory of its current state.  At each time step, each agent updates
its state according to a convex combination of its current state and
the information received from its neighbors. Critically, the weight
accorded to its own state tends to $1$ as time grows (as a way to
model the accumulation of experience).  The earliest SA algorithms
were proposed by Robbins and Monro \cite{VB-SM:51} who aimed to solve
root finding problems.  SA algorithms have then attracted much
interest due to many applications such as the study of reinforcement
learning~\cite{JNT:94}, consensus protocols in multi-agent
systems~\cite{GL-CG:17}, and fictitious play in game
theory~\cite{JH-WHS:02}.  A main tool in the study of SA algorithms
(see~\cite[Chapter 5]{HJK-GGY:97}) is the ordinary differential
equations (ODE) method, which transforms the analysis of asymptotic
properties of a discrete-time stochastic process into the analysis of
a continuous-time deterministic process.

In this paper, we consider linear SA algorithms
with random linear operators; these models are basic
first-order protocols with numerous applications in engineering and
sociology.
Currently, there are two main threads on the
 theoretical research of linear SA algorithms. One thread is based on assumptions that guarantee the state of the system
 converges to a deterministic point  \cite{HFC:96,MAK:96,EC-IW-SK:99,VBT:04,MAK-SS:15}. Another thread is the research on consensus of multi-agent systems, where the system matrices are assumed to be row-stochastic \cite{TL-JFZ:10,MH:12,GL-CG:17}.
 These two threads only consider a part of linear operators, and the critical condition for convergence is still unknown.
This paper
 develops appropriate analysis methods for  linear SA algorithms and also provides some sufficient and necessary conditions for convergence which
 include critical conditions for convergence of linear operators.
It is shown that under critical convergence conditions  the state of
the system will converge to random vectors, which is applied to consensus algorithms
over signed networks.
Moreover,  an additional
restriction of traditional SA algorithms is that only non-negative
gain schedules are allowed. This paper relaxes this
requirement and provides necessary and sufficient conditions for
convergence of linear SA algorithms under
arbitrary gains. In addition, we analyze the convergence rate of the
system when it is convergent.

Our general theoretical results are directly applicable to certain
multi-agent systems. The first application is to the study of
consensus problems in multi-agent systems. As it is well known,
numerous works provide sufficient conditions for consensus in
time-varying multi-agent systems with row-stochastic interaction
matrices; an incomplete list of references is
\cite{LM:05,FF-SZ:08a,ATS-AJ:08,TL-JFZ:10,GC-ZL-LG:14,GL-CG:17}; see
also the classic works~\cite{SC-ES:77,RC:84,JNT-DPB-MA:86}.
Recently, motivated by the study of antagonistic interactions in
social networks, novel concepts of bipartite, group, and cluster
consensus have been studied over signed networks (mainly focusing on
continuous-time dynamical models); see
\cite{CA:13,JY-LW:10,JQ-CY:13,ZG-KMY-KHJ-MC-YH:16}. In this paper, we
apply and extend our results on linear SA algorithms to the setting of
first-order discrete-time multi-agent system over random signed
networks and with state-dependent noise; for such models, we provide
novel necessary and sufficient conditions to reach consensus and group
consensus.

As the second application of our results, we study the
Friedkin-Johnsen (FJ) model of opinion dynamics in social
networks. The FJ model was first proposed in \cite{NEF-ECJ:99}, where
each agent is assumed to be susceptible to other agents' opinions but
also to be anchored to his own initial opinion with a certain level of
stubbornness. Ravazzi \emph{et al.} proposed a gossip version of the
FJ model in \cite{CR-PF-RT-HI:15}, whereby each link in the network is
sampled uniformly and the agents associated with the link meet and
update their opinions. The agents' opinions were proven to converge in
mean square. Frasca \emph{et al.} considered a symmetric pairwise
randomization of FJ in \cite{PF-HI-CR-RT:15}, whereby a pair of agents
are chosen to update their opinions. Our work, by exploiting
stochastic approximation, largely relaxes the conditions for
convergence when applied to FJ model over random interaction
networks. The sociological meaning of stochastic approximated FJ model
is that agents have cumulative memory about their previous
opinions. The adoption of SA models in the study of human behavior is
widely adopted in game theory and economics; e.g.,
see~\cite{JH-WHS:02}.

The main contributions of this paper are summarized as follows.
\begin{enumerate}
\item
For linear SA systems, we provide some necessary and sufficient conditions  to
  guarantee convergence by developing appropriate methods different from previous works.
  We derive some critical convergence conditions for linear operators for the first time.
    The convergence rate is also obtained when the system is convergent. Moreover, we consider the convergence of linear SA systems whose gain functions can take arbitrary real numbers.

\item Using our results, we get the necessary and sufficient
  conditions to reach consensus and group consensus of the first-order
  discrete-time multi-agent system over random signed networks and
  with state-dependent noise for the first time.

\item We extend our results to the multi-dimensional linear SA algorithms and provide applications to the multi-dimensional FJ model
  over random interaction networks.
\end{enumerate}

\paragraph*{Organization}
The remainder of this paper is organized as follows. We briefly review
the time-varying linear dynamical systems and propose a stochastic
approximation version of it in Section \ref{Protocol1}. The main
results are presented in Section \ref{Main_results}. In particular, we
introduce some preliminaries and assumptions in Subsection
\ref{sec:assumption}. Sufficient conditions that guarantee the
convergence of linear SA algorithms are obtained in Subsection
\ref{sf_2}. We provide the results on convergence rate in the same
subsection. In Subsection \ref{Necessary_Conditions}, we prove that
the sufficient condition is also necessary.  The necessary and
sufficient conditions for convergence are then summarized in
Subsection \ref{NScondition}. We generalize the results to
multi-dimensional models and discuss their application to group
consensus and the FJ model in Section \ref{sec:application}.  Section
\ref{sec:conclusion} concludes the paper.


\renewcommand{\thesection}{\Roman{section}}
\section{Linear Dynamical Systems}\label{Protocol1}
\renewcommand{\thesection}{\arabic{section}}

\subsection{Review of a time-varying linear dynamical system}

In \cite{CR-PF-RT-HI:15,YH-WL-TC:13} a time-varying linear dynamical system was
considered as follows:
\begin{eqnarray}\label{AF}
x(s+1)=P(s) x(s)+u(s),~~s=0,1,\ldots,
\end{eqnarray}
where $P(s)\in\real^{n\times n}$ is a matrix associated to the
communication network between agents, and $u(s)\in\real^{n}$ is an
input vector. Given a matrix $A\in\real^{n\times n}$, let $\rho(A)$
denote its spectral radius, i.e., $\rho(A)=\max_i |\lambda_i(A)|$,
where $\lambda_i(A)$ is an eigenvalue of $A$.  For system (\ref{AF}),
if $P(s)\equiv P$, $u(s) \equiv u$, and $\rho(A)<1$, then it is
immediate to see that $x(s)$ converges to $(I_{n}-P)^{-1} u$.

In this paper we will consider the case when  $\{P(s)\}$ and $\{u(s)\}$ are stochastic matrices and vectors respectively.
We define the $\sigma$-algebra generated by $\{P(s)\}$ and $\{u(s)\}$ as
 $\mathcal{F}_t=\sigma((P(s),u(s)),0\leq s\leq t).$ The probability space is
$(\Omega,\mathcal{F}_{\infty},P)$.

Since the system (\ref{AF}) does not
necessarily converge when $\{P(s)\}$ and $\{u(s)\}$ are stochastic, as an alternative, Ravazzi \emph{et
  al.}~\cite{CR-PF-RT-HI:15} investigate the ergodicity of system
(\ref{AF}) as follows.
\begin{proposition}[Theorem 1 in \cite{CR-PF-RT-HI:15}]\label{Ravazzi_theorem}
 Consider system (\ref{AF}) and assume $\{P(s)\}$ and $\{u(s)\}$ are
 sequences of
independent identically distributed (i.i.d.) random matrices and vectors with finite first moments. Assume there exists a
  constant $\alpha\in(0,1]$, a matrix $P\in\real^{n\times n}$ and a vector $u\in\real^{n}$  such that
    \begin{eqnarray*}
      \E[P(s)]=(1-\alpha)I_{n}+\alpha P,~~\E[u(s)]=\alpha u,~~\forall s\geq 0.
    \end{eqnarray*}
    If $\rho(P)<1$, then $x(s)$ converges to a random variable in
    distribution, and $\frac{1}{s} \sum_{k=0}^{s-1} x(k)$ converges to
    $(I_{n}-P)^{-1} u$ almost surely.
\end{proposition}


In this paper we adopt the stochastic approximation method to average
the effect of the stochastic $P(s)$ and $u(s)$ to the state $x(s)$.
In this case we study the sufficient and necessary conditions for
convergence of $x(s)$, and also obtain a convergence rate.

\subsection{Linear SA algorithms over random networks}
In this subsection we consider the stochastic-approximation version of
system (\ref{AF}), formulated as:
\begin{multline}\label{SAF_m1}
x(s+1)=(1-a(s))x(s) \\ +a(s) [P(s) x(s)+u(s)], \quad s=0,1,\ldots,
\end{multline}
where $a(s)\in\real$ is the gain function. The system (\ref{AF}) is so
called as linear SA algorithms
\cite{HFC:96,EC-IW-SK:99,VBT:04,TL-JFZ:10,MH:12,GL-CG:17}.  Compared
to system (\ref{AF}), each agent in system (\ref{AF}) updates its
state depending not only on the linear map $P(s) x(s)+u(s)$ but also
on its own current state.  If $a(s)=\frac{1}{s+1}$, then $x(s+1)$
equals the approximate average value of the previous $s$ linear maps
because $x(s)$ carries the information of the previous $s-1$ linear
maps.  Intuitively, in this case $x(s)$ approximately
  equals $\frac{1}{s} \sum_{k=0}^{s-1} x(k)$ in system (\ref{AF}), so
  that it should have the same limit as in Proposition
  \ref{Ravazzi_theorem}.  In fact, this result can be deduced by the
  following Proposition \ref{SAF_sc}.  Of course, this paper considers
  the more general case of $\{a(s)\}$ and $\{P(s)\}$.

The system (\ref{SAF_m1}) is a basic first-order discrete-time
multi-agent system with much prior theoretical analysis.
 A main thread in the research of such a system is to
  study the setting in which $x(s)$ converges to a deterministic
  point.  In \cite{HFC:96,EC-IW-SK:99}, convergence and convergence
  rates are studied for bounded linear operators with the assumption
  that there exists a matrix $P\in\real^{n\times n}$ whose
  eigenvalues' real parts are all less than $1$ such that
\begin{eqnarray}\label{HFC_con}
\lim_{s\rightarrow\infty}\Big(\sup_{s\leq t \leq m(s,T)} \Big\|\sum_{i=s}^t a(i)(P(i)-P) \Big\|_2\Big)=0,
\end{eqnarray}
where $m(s,T):=\max\{k:a(s)+\cdots+a(k)\leq T\}$ with $T$ being an
arbitrary positive constant, and $\|\cdot\|_2$ denotes the Euclidean
norm.  Later, Tadi\'{c} relaxed the boundary condition of $P(s)$ and
provided some convergence rates based on (\ref{HFC_con}) and the
assumption that the real parts of the eigenvalues of $P+\alpha I_n$
are all less than $1$, where $\alpha$ is a positive constant
\cite{VBT:04}.  Additionally, there are results on convergence rates
by assuming that $\{I_n-P(s)\}_{s\geq 0}$ are a sequence of positive
semi-definite matrices and $I_n-P$ is a positive definite matrix
\cite{MAK:96,MAK-SS:15}.  Another thread in the theoretical research
on system (\ref{SAF_m1}) is to consider its consensus behavior where
$\{P(s)\}$ and $\{u(s)\}$ are assumed to be row-stochastic matrices
and zero-mean noises respectively \cite{TL-JFZ:10,MH:12,GL-CG:17}. In
addition, system~(\ref{SAF_m1}) has many applications like computation of
PageRank \cite{WXZ-HFC-HTF:13}, sensor localization of wireless
networks \cite{UAK-SK-JMFM:09}, distributed consensus of multi-agent
systems, and belief evolution on social networks.

Despite all this prior theoretical research on system (\ref{SAF_m1}),
a key problem remains unsolved: What is the necessary and sufficient
condition for convergence regarding $\{P(s)\}$ and $u(s)$? Previous
works focused on the case when the real parts of the eigenvalues of
$P$ are all assumed to be less than $1$
\cite{HFC:96,EC-IW-SK:99,VBT:04,TL-JFZ:10,MH:12,GL-CG:17}, but it is
not known what happens when this condition is not satisfied. Also,
traditional SA algorithms consider only non-negative gains, so another
interesting problem is to investigate what happens if the gain
function $a(s)$ can take arbitrary real numbers. This paper considers
these two problems and studies the mean-square convergence of $x(s)$,
whose definition is given as follows:
\begin{definition}\label{MS_def}
  For an $n$-dimensional random vector $x$, we say $x(s)$ converges to
  $x$ in mean square if
  \begin{eqnarray}\label{MS_def_1}
    \E\|x\|_2^2<\infty~~\mbox{and}~~\lim_{s\to\infty}\E\|x(s)-x\|_2^2=0.
  \end{eqnarray}
  Also, we say $\{x(s)\}$ is mean-square convergent if there exists an
  $n$-dimensional random vector $x$ such that (\ref{MS_def_1}) holds.
\end{definition}

\renewcommand{\thesection}{\Roman{section}}
\section{Main results}\label{Main_results}
\renewcommand{\thesection}{\arabic{section}}

\renewcommand{\thesection}{\Roman{section}}
\subsection{Informal statement of main results}\label{sec:assumption}
\renewcommand{\thesection}{\arabic{section}}
We start with some notation.  Given a matrix $A\in\real^{n\times
  n}$, define
$\widetilde{\rho}_{\max}(A):=\max_{i}\mbox{Re}(\lambda_i(A))$ and
$\widetilde{\rho}_{\min}(A):=\min_{i}\mbox{Re}(\lambda_i(A))$ to be
the maximum and minimum values of the real parts of the eigenvalues of
$A$ respectively. It is easy to show that
$|\widetilde{\rho}_{\max}(A)|\leq \rho(A)$.

For $\{P(s)\}$ and $\{u(s)\}$, we relax the i.i.d. condition in
\cite{CR-PF-RT-HI:15} to the following assumption:

\textbf{(A1)} Suppose there exist a matrix $P\in\real^{n\times
  n}$ and a vector $u\in\real^n$ such that $\E[P(s)\,|\,x(s)]=P$ and
$\E[u(s)\,|\,x(s)]=u$ for any $s\geq 0$ and $x(s)\in\real^n$. Also,
assume $\E[\|P(s)\|_2^2\,|\,x(s)]$ and $\E[\|u(s)\|_2^2\,|\,x(s)]$ are uniformly
bounded.

For $\{a(s)\}$, generally SA algorithms use the following assumption:

\textbf{(A2)} Assume $\{a(s)\}$ are non-negative real numbers independent with $\{x(s)\}$, and satisfying  $\sum_{s=0}^{\infty} a(s)=\infty$ and $\sum_{s=0}^{\infty} a^2(s)<\infty$.

We will also consider the following alternative assumption.

\textbf{(A2')}  Assume $\{a(s)\}$ are non-positive real numbers independent with $\{x(s)\}$, and satisfying  $\sum_{s=0}^{\infty} a(s)=-\infty$ and $\sum_{s=0}^{\infty} a^2(s)<\infty$.

Under the assumptions (A1) and (A2), the previous works has investigated the cases when $\widetilde{\rho}_{\max}(P)<1$ and  $P(s)x+u(s)$ is a bounded linear operator for all
$s\geq 0$  \cite{HFC:96,EC-IW-SK:99}, or
 $\widetilde{\rho}_{\max}(P+\alpha I_n)<1$ \cite{VBT:04}, or $\{P(s)\}$ are row-stochastic matrices and $u=\bf{0}$ \cite{TL-JFZ:10,MH:12,GL-CG:17}. This paper will consider all the cases of $P$ and $u$, and show
 the necessary and
sufficient condition for the convergence of $x(s)$ in system
(\ref{SAF_m1}) is $\widetilde{\rho}_{\max}(P)<1$, or
$\widetilde{\rho}_{\max}(P)=1$ together with the following
condition for $P$ and $u$:

 \textbf{(A3)} Assume any eigenvalue of $P$ whose real part is $1$ equals $1$,
and the eigenvalue $1$ has the same algebraic and geometric multiplicities, and $\xi^T u=0$ for any left eigenvector $\xi^T$ of $P$ corresponding to the eigenvalue $1$.

Similarly, under (A1) and (A2') the necessary and sufficient condition
for the convergence of $x(s)$ is $\widetilde{\rho}_{\min}(P)>1$, or
$\widetilde{\rho}_{\min}(P)=1$ with (A3).

Also, we will study the convergence rates when $x(s)$ is convergent,
and the convergence conditions when $\{a(s)\}$ are arbitrary real
numbers.

\renewcommand{\thesection}{\Roman{section}}
\subsection{Sufficient convergence conditions and convergence rates}\label{sf_2}
\renewcommand{\thesection}{\arabic{section}}

Recall that $P$ and $u$ are the expectations of $P(s)$ and $u(s)$ respectively.
Let
\begin{eqnarray}\label{Jordan}
P=H^{-1}\mbox{diag}(J_1,\ldots,J_K)H:=H^{-1} D H,
\end{eqnarray}
 where $H\in\complex^{n\times n}$ is an invertible matrix, and $D$ is the Jordan normal form of $P$
with
\begin{eqnarray*}\label{Jordan_1}
J_i=
\begin{bmatrix}
    \lambda_{i'}(P) & 1 &  &     \\
     & \lambda_{i'}(P) & \ddots &     \\
     &  & \ddots & 1   \\
     &  &  &   \lambda_{i'}(P)
\end{bmatrix}_{m_i\times m_i}
\end{eqnarray*}
for $1\leq i\leq K$,  where $\lambda_{i'}(P)$ is the eigenvalue of $P$ corresponding to the Jordan block $J_i$.

Let $r$ be the algebraic multiplicity of the eigenvalue $1$ of $P$.
We first consider the case $\widetilde{\rho}_{\max}(P)=1$ (or
$\widetilde{\rho}_{\min}(P)=1$) with (A3), which implies that $r\geq
1$ and that the geometric multiplicity of the eigenvalue $1$ is
equal to $r$.  We choose a suitable $H$ such that
$\lambda_1(P)=\cdots=\lambda_r(P)=1$. Then the Jordan normal form $D$
can be written as
\begin{eqnarray}\label{Jordan_temp}
  D=
  \begin{bmatrix}
    I_r & {\mathbf{0}}_{r\times(n-r)}     \\
    {\mathbf{0}}_{(n-r)\times r} & \underline{D}_{(n-r)\times(n-r)}  \\
  \end{bmatrix}\in\complex^{n\times n},
\end{eqnarray}
where $\underline{D}:=\mbox{diag}(J_{r+1},\ldots,J_K)\in
\complex^{(n-r)\times (n-r)}$.  For any vector $y\in\complex^n$,
throughout this subsection we set $\bar{y}:=(y_1,\ldots,y_r)^\top$ and
$\underline{y}:=(y_{r+1},\ldots,y_n)^\top$.

\begin{theorem}\label{SAF_c1}(Convergence of linear SA algorithms at critical point)
  Consider the system (\ref{SAF_m1}) satisfying (A1), (A2), and (A3)
  with $\widetilde{\rho}_{\max}(P)=1$, or satisfying (A1), (A2'), and
  (A3) with $\widetilde{\rho}_{\min}(P)=1$.  Let $H$ be the matrix
  defined by (\ref{Jordan}) such that the Jordan normal form $D$ has
  the form of (\ref{Jordan_temp}).  Then, for any initial state,
  $x(s)$ converges to $H^{-1} y$ in mean square, where $\bar{y}$ is a
  random vector satisfying $\E\bar{y}=\overline{Hx(0)}$ and
  $\E\|\bar{y}\|_2^2<\infty$, and
  $\underline{y}=(I_{n-r}-\underline{D})^{-1} \underline{Hu}$.
\end{theorem}

From Theorem \ref{SAF_c1}, $x(s)$ converges to a random vector under the critical condition  $\widetilde{\rho}_{\max}(P)=1$ (or $\widetilde{\rho}_{\min}(P)=1$), which is different from the previous works where $x(s)$ converges to a deterministic vector under non critical conditions \cite{HFC:96,EC-IW-SK:99,VBT:04,TL-JFZ:10,MH:12,GL-CG:17}.
Due to this difference, the traditional method cannot be used in the proof of Theorem \ref{SAF_c1}. We propose a new method to prove this theorem as follows.

\begin{IEEEproof}[Proof of Theorem \ref{SAF_c1}]
Let $y(s):=H x(s)$, $v(s):=H u(s)$ and $D(s):=H P(s) H^{-1}$, then by (\ref{SAF_m1}) we have
\begin{multline}\label{SAF_c1_4}
H^{-1} y(s+1)\\=(1-a(s))H^{-1} y(s)+a(s)[P(s) H^{-1} y(s)+ u(s)],
\end{multline}
which implies
\begin{eqnarray}\label{SAF_c1_5}
y(s+1)=y(s)+a(s)[(D(s)-I_n)y(s)+v(s)].
\end{eqnarray}
Let $v:=\E[v(s)]=Hu$. From (\ref{Jordan}) we have $HP=DH$, which implies $H_i P= H_i$ for $1\leq i\leq r$, where $H_i$ is the $i$-th row of the matrix $H$.
Thus, $H_i$, $1\leq i\leq r$, is a left eigenvector corresponding to the eigenvalue $1$. By (A3) we have
\begin{eqnarray}\label{SAF_c1_5_1}
v_i=H_i u=0, ~~~~\forall 1\leq i\leq r.
 \end{eqnarray}
Recall that $\underline{v}=(v_{r+1},\ldots,v_n)^\top$. Also, $I_{n-r}-\underline{D}$ is an invertible matrix, so we can set
\begin{eqnarray*}\label{SAF_c1_6}
z:=
\begin{bmatrix}
     \mathbf{0}_{r\times 1}     \\
     (I_{n-r}-\underline{D})^{-1}\underline{v}  \\
\end{bmatrix}\in \complex^n.
\end{eqnarray*}
From (\ref{Jordan_temp}) and (\ref{SAF_c1_5_1}) we have
\begin{eqnarray}\label{SAF_c1_6_1}
(D-I_n)z+v=\mathbf{0}_{n\times 1}.
\end{eqnarray}

Set $\theta(s):=y(s)-z$. From (\ref{SAF_c1_5}) we obtain
\begin{eqnarray}\label{SAF_c1_7}
\theta(s+1)=\theta(s)+a(s)[(D(s)-I_n)(\theta(s)+z)+v(s)].
\end{eqnarray}

We first consider the case when $\widetilde{\rho}_{\max}(P)=1$, which
implies that $\underline{D}-I_{n-r}$ is a Hurwitz matrix.  Thus, by the
stability theory of continuous Lyapunov equation
(see~\cite[Corollary~2.2.4]{RAH-CRJ:94}), there exists a Hermitian
positive definite matrix $A\in\complex^{(n-r)\times (n-r)}$ such that
\begin{eqnarray}\label{SAF_c1_8}
\begin{aligned}
(\underline{D}-I_{n-r}) ^* A+ A(\underline{D}-I_{n-r})=-I_{n-r},
\end{aligned}
\end{eqnarray}
where $(\cdot)^*$ denotes the  conjugate transpose of the matrix or vector.
Set
\begin{eqnarray*}\label{SAF_c1_9}
A_1:=
\begin{bmatrix}
    I_r & {\mathbf{0}}_{r\times(n-r)}     \\
     {\mathbf{0}}_{(n-r)\times r} & A_{(n-r)\times(n-r)} \\
\end{bmatrix}\in \complex^{n\times n},
\end{eqnarray*}
then $A_1$ is still a Hermitian positive definite matrix.
Define the  Lyapunov function $V_1(\theta):=\theta^* A_1 \theta$.
By (\ref{SAF_c1_7}), (A1) and (\ref{SAF_c1_8}), for any $\theta(s)$ we have
\begin{align}\label{SAF_c1_10}
\E [ &V_1(\theta(s+1))|\theta(s)]  \\
&\leq V_1(\theta(s))+a(s)\theta^* (s)\big[ (D-I_n)^* A_1+A_1 (D-I_n) \big] \theta(s)\nonumber \\
&\quad+O\big(a^2(s)(\|\theta(s)\|_2^2+1)\big)\footnotemark.\nonumber
\end{align}
 \footnotetext{Given two sequences of positive numbers $\{g_1(s)\}_{s=0}^{\infty}$ and
$\{g_2(s)\}_{s=0}^{\infty}$, we say $g_1(s) = O(g_2(s))$ if there exist a
constants $c>0$ such that $g_1(s)\leq c g_2(s)$ for
all $s\geq 0$.}
From (\ref{Jordan_temp}) and (\ref{SAF_c1_8}), we obtain
\begin{align}
(&D-I_n)^* A_1 + A_1 (D-I_n)\nonumber \\
&=\begin{bmatrix}
    {\mathbf{0}}_{r\times r} & {\mathbf{0}}_{r\times (n-r)}     \\
     {\mathbf{0}}_{(n-r)\times r} & (\underline{D}-I_{n-r})^* A+A(\underline{D}-I_{n-r}) \\
\end{bmatrix}\nonumber\\
&=\begin{bmatrix}
    {\mathbf{0}}_{r\times r} & {\mathbf{0}}_{r\times (n-r)}     \\
     {\mathbf{0}}_{(n-r)\times r} & -I_{n-r} \\
\end{bmatrix},\label{SAF_c1_11}
\end{align}
so (\ref{SAF_c1_10}) implies
\begin{eqnarray}\label{SAF_c1_12}
\E[V_1(\theta(s+1))]\leq [1+c_1 a^2(s)] \E[V_1(\theta(s))]+c_2 a^2(s),
\end{eqnarray}
where $c_1$ and $c_2$ are two positive constants.
Using (\ref{SAF_c1_12}) repeatedly we get
\begin{align}\label{SAF_cl_13}
\E[&V_1(\theta(s+1))]\\
&\leq \prod_{i=0}^s [1+c_1 a^2(i)] +\sum_{i=0}^s c_2 a^2(i)\prod_{j=i+1}^s [1+c_1 a^2(j)]\nonumber\\
&<\infty ~~~~\mbox{as}~~s\to\infty,\nonumber
\end{align}
where the last inequality uses the condition that $\sum_{s=0}^{\infty} a^2(s)<\infty$.
Also, because $A_1$ is a Hermitian positive definite matrix,
\begin{eqnarray}\label{SAF_c1_13_1}
\frac{1}{\rho(A_1)} V_1(\theta(s))\leq \|\theta(s)\|_2^2\leq \frac{1}{\lambda_{\min}(A_1)}V_1(\theta(s)).
\end{eqnarray}
Combining (\ref{SAF_cl_13}) and
(\ref{SAF_c1_13_1}) yields
\begin{eqnarray}\label{SAF_cl_14}
\sup_{s} \E\|\theta(s)\|_2^2\leq \sup_s \frac{\E[V_1(\theta(s))]}{\lambda_{\min}(A_1)}<\infty.
\end{eqnarray}

Inequality (\ref{SAF_cl_14}) shows that $\theta(s)$ will not diverge, however we need to prove its convergence.
We first consider the convergence of $\underline{\theta}(s)$.
Set
\begin{eqnarray*}\label{SAF_c1_15}
A_2:=
\begin{bmatrix}
    {\mathbf{0}}_{r\times r} & {\mathbf{0}}_{r\times (n-r)}     \\
     {\mathbf{0}}_{(n-r)\times r} & A_{(n-r)\times(n-r)} \\
\end{bmatrix}\in \complex^{n\times n}
\end{eqnarray*}
and define $V_2(\theta):=\theta^* A_2 \theta=\underline{\theta}^* A
\underline{\theta}$. Similar to (\ref{SAF_c1_10}), we have
\begin{align}\label{SAF_c1_16}
\E[&V_2(\theta(s+1))]\\
&\leq \E\Big[V_2(\theta(s))+a(s)\theta^*(s) (D-I_n)^* A_2 \theta(s)\nonumber \\
&\quad+a(s) \theta^*(s) A_2 (D-I_n)\theta(s)+O\big(a^2(s)(\|\theta(s)\|_2^2+1)\big)\Big]\nonumber\\
&=\E\Big[V_2(\theta(s))+a(s)\theta^*(s)
\begin{bmatrix}
    {\mathbf{0}}_{r\times r} & {\mathbf{0}}_{r\times (n-r)}     \\
     {\mathbf{0}}_{(n-r)\times r} & -I_{n-r} \\
\end{bmatrix}
\theta(s)\Big]\nonumber\\
&\quad+O(a^2(s))   \nonumber\\
&\leq \Big(1-\frac{a(s)}{\rho(A)}\Big) \E[V_2(\theta(s))] +O\big(a^2(s)\big),\nonumber
\end{align}
where the forth line uses (\ref{SAF_c1_11}) and (\ref{SAF_cl_14}), and
the last inequality does a similar computation as (\ref{SAF_c1_13_1}).
By (\ref{SAF_c1_16}) and Lemma \ref{lguo} in Appendix
\ref{App_lemmas}, we obtain $\lim_{s\to\infty} \E[V_2(\theta(s))]=0$,
which implies
\begin{eqnarray}\label{SAF_c1_17}
\lim_{s\to\infty}\E\|\underline{\theta}(s)\|_2^2=0.
\end{eqnarray}

It remains to consider the convergence of $\bar{\theta}(s)$.
Set
\begin{eqnarray*}\label{SAF_c1_18}
A_3:=
\begin{bmatrix}
    I_r & {\mathbf{0}}_{r\times(n-r)}     \\
     {\mathbf{0}}_{(n-r)\times r} & {\mathbf{0}}_{(n-r)\times(n-r)} \\
\end{bmatrix}\in \real^{n\times n},
\end{eqnarray*}
and define $V_3(\theta)=\theta^* A_3 \theta= \bar{\theta}^* \bar{\theta}$. By (\ref{Jordan_temp}) and (\ref{SAF_c1_5_1}) we get
$A_3 (D-I_n)=\mathbf{0}_{n\times n}$ and $A_3 v=\mathbf{0}_{n\times 1}$, thus by (A1) for any $i<j$ we have
\begin{eqnarray}\label{SAF_c1_19}
\begin{aligned}
&E\Big[[(D(i)-I_n)(\theta(i)+z)+v(i)]^* A_3 \\
&\qquad\times [(D(j)-I_n)(\theta(j)+z)+v(j)]\Big]\\
&=E\Big[[(D(i)-I_n)(\theta(i)+z)+v(i)]^* A_3 \\
&\qquad\times [(D-I_n)(\theta(j)+z)+v]\Big]=0.
\end{aligned}
\end{eqnarray}
Similarly, the equation (\ref{SAF_c1_19}) still holds for $i>j$.
From these and (\ref{SAF_c1_7}) we get for any $s_2> s\geq 0$,
\begin{align}\label{SAF_c1_20}
\E[V_3(\theta&(s_2)-\theta(s))]  \\
&=\E\Big[V_3\Big(\sum_{i=s}^{s_2-1}[\theta(i+1)-\theta(i)]\Big)\Big]\nonumber\\
&=\E\bigg[ \Big(\sum_{i=s}^{s_2-1} a(i)[(D(i)-I_n)(\theta(i)+z)+v(i)]\Big)^* A_3\nonumber \\
&\quad\times  \Big(\sum_{i=s}^{s_2-1} a(i)[(D(i)-I_n)(\theta(i)+z)+v(i)]\Big) \bigg]\nonumber\\
&=\sum_{i=s}^{s_2-1}  a^2(i) \E\big[ [(D(i)-I_n)(\theta(i)+z)+v(i)]^*  \nonumber\\
&\quad\times A_3 [(D(i)-I_n)(\theta(i)+z)+v(i)] \big]\nonumber\\
&= O\Big(\sum_{i=s}^{s_2-1}a^2(i)\Big),\nonumber
\end{align}
where the last line uses (A1) and (\ref{SAF_cl_14}). Since $\sum_{i=0}^{\infty} a^2(i)<\infty$, from
(\ref{SAF_c1_20}) we have
\begin{multline}\label{SAF_c1_21}
\lim_{s\to\infty}\lim_{s_2\to\infty} \E\|\bar{\theta}(s_2)-\bar{\theta}(s)\|_2^2\\
=\lim_{s\to\infty} \lim_{s_2\to\infty} \E[V_3(\theta(s_2)-\theta(s))]=0.
\end{multline}
By the Cauchy criterion (see~\cite[page~58]{JAH:07}),
$\bar{\theta}(s)$ has a mean square limit $\bar{\theta}(\infty)$.
Also, from (\ref{SAF_c1_7}), (A1) and (\ref{SAF_c1_5_1}) we have
\begin{align}\label{SAF_c1_22}
\E[A_3 \theta&(s+1)]\\
&=\E\big[\E[A_3 \theta(s+1)\,|\,\theta(s)]\big]\nonumber\\
&=\E\big[A_3 \theta(s)+a(s)A_3[(D-I_n)(\theta(s)+z)+v]\big]\nonumber\\
&=\E[A_3 \theta(s)]=\cdots= A_3 \theta(0),\nonumber
\end{align}
which is followed by
\begin{eqnarray}\label{SAF_c1_23}
\E\bar{\theta}(\infty)=\bar{\theta}(0)=\bar{y}(0)=\overline{Hx}(0).
\end{eqnarray}

We remark that $x(s)=H^{-1}[\theta(s)+z]$. Let $y$ be a vector
satisfying
$\underline{y}=\underline{z}=(I_{n-r}-\underline{D})^{-1}\underline{v}$ and
$\bar{y}=\bar{\theta}(\infty)+\bar{z}=\bar{\theta}(\infty).$ By
(\ref{SAF_c1_17}) and (\ref{SAF_c1_21}) we have that $x(s)$ converges to
$H^{-1}y$ in mean square. By (\ref{SAF_c1_23}) and (\ref{SAF_cl_14})
we get $\E \bar{y}=\overline{Hx}(0)$ and $\E\|\bar{y}\|_2^2<\infty$.

For the case that $\widetilde{\rho}_{\min}(P)=1$, which implies
$I_{n-r}-\underline{D}$ is a Hurwitz matrix. Set $b(s)=-a(s)\geq 0$ and substitute it to (\ref{SAF_c1_7}) we obtain
\begin{eqnarray*}\label{SAF_c1_24}
\theta(s+1)=\theta(s)+b(s)[(I_n-D(s))(\theta(s)+z)-v(s)].
\end{eqnarray*}
Finally, a process similar to that from (\ref{SAF_c1_8}) to
(\ref{SAF_c1_23}) yields our result.
\end{IEEEproof}

For the case when $\widetilde{\rho}_{\max}(P)<1$ or $\widetilde{\rho}_{\min}(P)>1$,
from the proof of Theorem \ref{SAF_c1}
we have the following proposition:
\begin{proposition}\label{SAF_sc}
Consider the system (\ref{SAF_m1}) satisfying (A1), (A2) and $\widetilde{\rho}_{\max}(P)<1$, or satisfying
(A1), (A2') and $\widetilde{\rho}_{\min}(P)>1$.
Then, for any initial state, $x(s)$ converges to $(I_n-P)^{-1} u$ in mean square.
\end{proposition}
\begin{IEEEproof}
We can set $r=0$ in the proof of Theorem \ref{SAF_c1}, then we obtain that $x(s)$ converges to $H^{-1}(I_n-D)^{-1}Hu=(I_n-P)^{-1}u$ in mean square.
\end{IEEEproof}

Next, we give the convergence rate when $x(s)$ is mean-square convergent.
\begin{theorem}\label{SAF_r1}(Convergence rates of linear SA algorithms)
Consider the system (\ref{SAF_m1}) satisfying (A1)
 and one of the following four cases: i) $\widetilde{\rho}_{\max}(P)<1$; ii) $\widetilde{\rho}_{\min}(P)>1$; iii) $\widetilde{\rho}_{\max}(P)=1$ with (A3);
and iv) $\widetilde{\rho}_{\min}(P)=1$ with (A3).
Let $\beta>0, \gamma\in (\frac{1}{2},1]$, and $\alpha$ be a large positive number.
Choose  $a(s)=\frac{\alpha}{(s+\beta)^{\gamma}}$ if $\widetilde{\rho}_{\max}(P)\leq 1$,
and $a(s)=\frac{-\alpha}{(s+\beta)^{\gamma}}$ if $\widetilde{\rho}_{\min}(P)\geq 1$.
Then for any initial state,
\begin{multline*} 
  \E\big\|x(s)-x\big\|_2^2\\
  = \begin{cases}
    O(s^{-\gamma}),  &\mbox{ if } \widetilde{\rho}_{\max}(P)<1 \mbox{ or } \widetilde{\rho}_{\min}(P)>1\\
    O(s^{1-2\gamma}), &\mbox{ if } \widetilde{\rho}_{\max}(P)=1 \mbox{ or } \widetilde{\rho}_{\min}(P)=1
  \end{cases}
\end{multline*}
where $x$ is a mean square limit of $x(s)$ whose expression is
provided by Theorem \ref{SAF_c1} and Proposition \ref{SAF_sc}.
\end{theorem}

The proof of this theorem is postponed to Appendix  \ref{proof_SAF_r1}.

\begin{remark}
 For the case when $\widetilde{\rho}_{\max}(P)<1$, there
  exist results on the convergence and convergence rates of $x(s)$
  provided some additional conditions hold, beside (A1)-(A2).  For
  example, if $\lim_{s\rightarrow\infty}\sum_{k=0}^{s-1}\frac{
    \|P(s)\|_2}{s}$ a.s. exists and $\widetilde{\rho}_{\max}(P+\alpha
  I_n)<1$ with $\alpha$ being a positive constant, then Theorem 2 in
  \cite{VBT:04} provides sufficient and necessary conditions for the
  convergence rate of $x(s)$; if $\|x(s)\|_2$ is uniformly bounded
  a.s., then by the ODE method in SA theory (Theorem 5.2.1 in
  \cite{HJK-GGY:97} or Theorem 2.2 in \cite{VSB-SPM:00}) we have
  $x(s)$ converges to $(I_n-P)^{-1} u$ a.s.  However, to the best of
  our knowledge, our results in Proposition \ref{SAF_sc} and Theorem
  \ref{SAF_r1} cannot be deduced from existing results without
  additional conditions.
\end{remark}
\renewcommand{\thesection}{\Roman{section}}
\subsection{Necessary conditions for convergence}\label{Necessary_Conditions}
\renewcommand{\thesection}{\arabic{section}}
We first consider necessary conditions of convergence under the assumptions (A1) and (A2) or (A2'):
\begin{theorem}\label{SAF_nca}
Consider the system (\ref{SAF_m1}) satisfying (A1). Then:\\
i) If $\widetilde{\rho}_{\max}(P)>1$, or
$\widetilde{\rho}_{\max}(P)=1$ but (A3) does not hold,
there exist some initial states such that $x(s)$ is not mean-square convergent for any
$\{a(s)\}$ satisfying (A2).\\
ii) If $\widetilde{\rho}_{\min}(P)<1$, or
$\widetilde{\rho}_{\min}(P)=1$ but (A3) does not hold,
there exist some initial states such that $x(s)$ is not mean-square convergent for any
$\{a(s)\}$ satisfying (A2').
\end{theorem}

The proof of this theorem is postponed to Appendix \ref{proof_nca}.

The necessary condition of convergence in Theorem \ref{SAF_nca} has a
constraint that the gain function $\{a(s)\}$ must satisfy the
assumption (A2) or (A2'). An interesting problem is to understand what
happens if $\{a(s)\}$ are chosen as arbitrary real numbers. Obviously,
from protocol (\ref{SAF_m1}) if $\{a(s)\}$ has only finite non-zero
elements, then $x(s)$ will converge to a random variable. Thus, we only
consider the setting whereby $x(s)$ does not converge to a
deterministic vector for arbitrary gains.

Recall that
\begin{eqnarray*}
P=H^{-1} \mbox{diag}(J_1,\ldots,J_K)H=H^{-1} D H,
\end{eqnarray*}
where $H\in\complex^{n\times n}$ is an invertible matrix, and $D$ is
the Jordan normal form of $P$.  For $1\leq i\leq K$, define
\begin{eqnarray}\label{Jordan_2}
\widetilde{I}_i=\mbox{diag}(0,\ldots, I_{m_i},\ldots,0)\in\real^{n\times n},
\end{eqnarray}
which corresponds to the Jordan block $J_i$ and then $D\widetilde{I}_i=\mbox{diag}(0,\ldots,J_i,\ldots,0)$.
To study the necessary condition
for convergence of  system (\ref{SAF_m1}), we need the following two assumptions:

\textbf{(A4)} Assume there is a Jordan block $J_{j}$ in $D$ associated with the eigenvalue $\lambda_{j'}(P)$
such that $\mbox{Re}(\lambda_{j'}(P))=1$ and
\begin{equation}\label{Assump_ag_1}
  \E\big[\| \widetilde{I}_j H[(P(s)-P)x(s)+u(s)-u]\|_2^2\,|\,x(s) \big]
  \geq c_1 \|x(s)\|_2^2+c_2
\end{equation}
for any  $s\geq 0$ and $x(s)\in\real^n$,
where $P$, $u$, $H$, $D$ and $\widetilde{I}_j$ are defined by (A1), (\ref{Jordan}), and (\ref{Jordan_2}), and
$c_1$ and $c_2$ are constants satisfying $c_1\geq 0$, $c_2\geq 0$, and $c_1+c_2>0$.

\textbf{(A4')} Assume there are two Jordan blocks $J_{j_1}$ and $J_{j_2}$ associated with the eigenvalues
$\lambda_{j_1'}(P)$ and $\lambda_{j_2'}(P)$ respectively
such that $\mbox{Re}(\lambda_{j_1'}(P))<1<\mbox{Re}(\lambda_{j_2'}(P))$ and
(\ref{Assump_ag_1}) holds for $j=j_1, j_2$.

\begin{theorem}\label{SAF_r2}
Consider the system (\ref{SAF_m1}) satisfying (A1) and (A4) or (A4'). In addition, assume there exists a constant $c_3>0$ such that
for any $s\geq 0$ and $x(s)\in\real^n$,
\begin{eqnarray}\label{SAF_r2_01}
\begin{aligned}
&\E\big[\|(P(s)-P)x(s)+u(s)-u\|_2^2|x(s)\big]\geq c_3.
\end{aligned}
\end{eqnarray}
Then for any deterministic vector $b\in\real^n$,  any initial state $x(0)\neq b$, and any real number sequence $\{a(s)\}_{s\geq 0}$ independent with
$\{x(s)\}_{s\geq 0}$, $x(s)$ cannot converge to $b$ in mean square.
\end{theorem}
The proof of this theorem is postponed to Appendix \ref{Proof_SAF_r2}.


If $u(s)$ is a degenerate random vector which means that $\E\|u(s)-u\|_2^2=0$, then the condition (\ref{SAF_r2_01}) may not be satisfied.

\begin{theorem}\label{SAF_r3}
Consider the system (\ref{SAF_m1}) satisfying (A1), and $\E[\|u(s)-u\|_2^2\,|\,x(s)]=0$
for any $s\geq 0$ and $x(s)\in\real^n$. Assume (A4) or (A4') holds but using
\begin{eqnarray}\label{SAF_r3_01}
\begin{aligned}
\E\big[\| \widetilde{I}_j H(P(s)-P)x(s)\|_2^2\,|\,x(s)\big]\geq c_1 \|x(s)\|_2^2
\end{aligned}
\end{eqnarray}
instead of (\ref{Assump_ag_1}).
For any deterministic vector $b\in\real^n$ and any initial state $x(0)\neq b$,
if  one of the following three conditions holds:\\
i) $u\neq {\mathbf{0}}_{n\times 1}$ and $x(0)\neq {\mathbf{0}}_{n\times 1}$;\\
ii) $u\neq {\mathbf{0}}_{n\times 1}$, $x(0)={\mathbf{0}}_{n\times 1}$, and $b \neq \alpha u $ for any $\alpha\in\real$; or\\
iii) $u={\mathbf{0}}_{n\times 1}$, and the eigenvalues $\lambda_{j'}(P)$ in (A4), or $\lambda_{j_1'}(P)$ and $\lambda_{j_2'}(P)$ in (A4') are not real numbers,\\
then  $x(s)$ cannot converge to $b$ in mean square for any real number sequence $\{a(s)\}_{s\geq 0}$ independent with
$\{x(s)\}_{s\geq 0}$.
\end{theorem}
The proof of this theorem is postponed to Appendix \ref{Proof_SAF_r3}.

\renewcommand{\thesection}{\Roman{section}}
\subsection{Necessary and sufficient conditions for convergence}\label{NScondition}
\renewcommand{\thesection}{\arabic{section}}

From  Theorems \ref{SAF_c1} and \ref{SAF_nca} and Proposition \ref{SAF_sc}, the following
necessary and sufficient condition for convergence with non-negative gains is obtained immediately.

\begin{theorem}\label{SAF_sn}(Necessary and sufficient condition for convergence of linear SA algorithms with non-negative gains)
Consider the system (\ref{SAF_m1}) satisfying (A1) and (A2).
Then $x(s)$ is mean-square convergent for any initial state
if and only if  $\widetilde{\rho}_{\max}(P)<1$, or $\widetilde{\rho}_{\max}(P)=1$ with (A3).
\end{theorem}

\begin{remark}
We remark that Theorem \ref{SAF_sn} is completely different from  previous sufficient and necessary conditions of convergence
in linear SA algorithms where only the case when $\widetilde{\rho}_{\max}(P)<1$ is considered and the assumptions are different from (A2) (Theorem 2 in \cite{HFC:96}; Theorem 1 in \cite{EC-IW-SK:99}; Theorems 1 and 2 in \cite{VBT:04}). In fact, the convergence of $x(s)$ at the critical point $\widetilde{\rho}_{\max}(P)=1$ has some applications such as the group consensus
over random signed networks; see Subsection \ref{subsec_gc}.
\end{remark}

Similarly, from  Theorems \ref{SAF_c1} and \ref{SAF_nca} and Proposition \ref{SAF_sc}, the following
necessary and sufficient condition for convergence with  non-positive gain is obtained immediately.

\begin{theorem}\label{SAF_sn2}(Necessary and sufficient condition for convergence of linear SA algorithms with non-positive gains)
Consider the system (\ref{SAF_m1}) satisfying (A1) and (A2').
Then $x(s)$ is mean-square convergent for any initial state
if and only if  $\widetilde{\rho}_{\min}(P)>1$, or $\widetilde{\rho}_{\min}(P)=1$ with (A3).
\end{theorem}


\begin{remark}
Compared to Theorem 1 in \cite{CR-PF-RT-HI:15}, Theorem \ref{SAF_sn}
extends the convergence condition from $\rho(P)<1$ to the sufficient
and necessary condition. In fact, for the basic linear dynamical system
$x(s+1)=P x(s)+u$, $x(s)$ converges if and only if $\rho(P)<1$.
However, if we consider the time-varying linear dynamical system and adopt the SA
method to eliminate the effect of fluctuation, then the convergence
condition can be substantially weakened.
\end{remark}

Theorems \ref{SAF_sn} and \ref{SAF_sn2} have a constraint that the gain
function $\{a(s)\}$ must satisfy the assumption (A2) or (A2').
Without this constraint we can get the following necessary and
sufficient condition for convergence to a deterministic vector, but
with some additional conditions on $\{u(s)\}$ or $\{P(s)\}$.

\begin{theorem}[Necessary and sufficient condition for convergence of linear SA algorithms with arbitrary gains]\label{SAF_r4}
Consider the system (\ref{SAF_m1}) which satisfies (A1).
Suppose there exists a constant $c\in(0,1)$ such that for any $s\geq 0$, $x(s)\in\real^n$,
 $\xi_1,\ldots,\xi_m \in \{P_{ij}(s),1\leq i,j\leq n; u_i(s), 1\leq i\leq n\}$
and $c_1,\ldots,c_m\in\complex$,
\begin{equation}\label{SAF_r4_01}
\E\Big[\Big|\sum_{i=1}^m c_i (\xi_i-\E \xi_i) \Big|^2\,|\,x(s)\Big]
\geq c\sum_{i=1}^m |c_i|^2 \E\big[(\xi_i-\E \xi_i)^2\,|\,x(s)\big].
\end{equation}
In addition, assume one of the following two conditions holds:\\
i) $\inf_{k,s} \E[(u_k(s)-u_k)^2\,|\,x(s)]>0$.\\
ii) $\E[\|u(s)-u\|^2\,|\,x(s)]=0$, $u\neq {\mathbf{0}}_{n\times 1}$, $x(0)\neq {\mathbf{0}}_{n\times 1}$,  and $\inf_{i,j,s} \E[(P_{ij}(s)-P_{ij})^2\,|\,x(s)]>0$.\\
Then we can choose a real number sequence $\{a(s)\}_{s\geq 0}$ independent with
$\{x(s)\}_{s\geq 0}$ such that $x(s)$ converges to a deterministic vector different from $x(0)$ in mean square
if and only if $\widetilde{\rho}_{\max}(P)<1$ or $\widetilde{\rho}_{\min}(P)>1$.
\end{theorem}
\begin{IEEEproof}
If $\widetilde{\rho}_{\max}(P)<1$ or $\widetilde{\rho}_{\min}(P)>1$, by  Proposition \ref{SAF_sc} we obtain that
$x(s)$ converges to $(I_n-P)^{-1}u$ in mean square.

For $\widetilde{\rho}_{\min}(P)\leq 1 \leq \widetilde{\rho}_{\max}(P)$, we
 set $\widetilde{P}(s):=P(s)-P$ and $\widetilde{u}(s):=u(s)-u$. Define $H$ and $K$ by
 (\ref{Jordan}), and define  $\widetilde{I}_i$ by (\ref{Jordan_2}).
 For any $j\in\{1,\ldots,K\}$,
since $H$ is an invertible matrix, $\widetilde{I}_j H$ contains at least one non-zero row $H_{j'}$.
Thus, for any
$x(s)\in\real^n$ we have
\begin{align}\label{SAF_r4_1}
\E\big[&\| \widetilde{I}_j H[\widetilde{P}(s)x(s)+\widetilde{u}(s)]\|_2^2\,|\,x(s)\big]\nonumber\\
&\geq \E\big[| H_{j'} [\widetilde{P}(s)x(s)+\widetilde{u}(s)]|^2\,|\,x(s)\big]\nonumber\\
&=  \E\Big[\Big|\sum_{i,k}  H_{j'i}  \widetilde{P}_{ik}(s)x_k(s)+\sum_{i} H_{j' i}\widetilde{u}_i(s)\Big|^2\,|\,x(s)\Big]\nonumber\\
&\geq c\sum_{i,k}  |H_{j'i}|^2 \E\big[\widetilde{P}_{ik}^2(s)\,|\,x(s)\big] x_k^2(s)\nonumber\\
&\quad +c\sum_{i} |H_{j' i}|^2 \E\big[\widetilde{u}_i^2(s)\,|\,x(s)\big],
\end{align}
where the last inequality uses (\ref{SAF_r4_01}).

If Condition i) holds, we have  there exists a constant $d_1>0$ such that
$\E\big[\widetilde{u}_i^2(s)\,|\,x(s)\big]\geq d_1$ for $s\geq 0$ and $1\leq i\leq n$.
 Combing this with (\ref{SAF_r4_1}) and the assumption $\widetilde{\rho}_{\min}(P)\leq 1 \leq \widetilde{\rho}_{\max}(P)$, we obtain that (\ref{SAF_r2_01}) and (A4) or (A4') hold.
By Theorem \ref{SAF_r2},  $x(s)$ cannot converge to a deterministic vector different from $x(0)$ in mean square.

If Condition ii) holds, we have $\E[\|\widetilde{u}(s)\|_2^2\,|\,x(s)]=0$ and
 there exists a constant $d_2>0$ such that
$\E\big[\widetilde{P}_{ik}^2(s)\,|\,x(s)\big]\geq d_2$ for $s\geq 0$ and $1\leq i,k\leq n$.
By (\ref{SAF_r4_1}) we obtain
\begin{align*}
\E\big[&\| \widetilde{I}_j H[\widetilde{P}(s)x(s)]\|_2^2\,|\,x(s)\big]\\
&\geq c d_2\sum_{i,k}  |H_{j'i}|^2 x_k^2(s)= c d_2 \|x(s)\|_2^2 \sum_{i}  |H_{j'i}|^2,
\end{align*}
which is followed by (\ref{SAF_r3_01}). By Theorem \ref{SAF_r3} i)  $x(s)$ cannot converge to a deterministic vector different from $x(0)$ in mean square.
\end{IEEEproof}

\renewcommand{\thesection}{\Roman{section}}
\section{Some Applications and Extension}\label{sec:application}
\renewcommand{\thesection}{\arabic{section}}

\renewcommand{\thesection}{\Roman{section}}
\subsection{Necessary and sufficient conditions for group consensus
  over random signed networks and with state-dependent noise}\label{subsec_gc}
\renewcommand{\thesection}{\arabic{section}}

As we discuss in the Introduction, consensus problems in multi-agent
systems have drawn a lot of attention from various fields including
physics, biology, engineering and mathematics in the past two decades.
Typically, a general assumption is adopted that the interaction matrix
associated with the network is row-stochastic at every time.  Recently,
motivated by the possible antagonistic interaction in social networks,
bipartite/group/cluster consensus problems have been studied over
signed networks (focusing on continuous-time dynamic models), e.g.,
see \cite{CA:13,JY-LW:10,JQ-CY:13,ZG-KMY-KHJ-MC-YH:16}.
On the other hand, SA has become a effective tool for the distributed consensus to
eliminate the effects of fluctuations \cite{MYH-JHM:09,RC-GC-PF-FG:11,NEL-AO:14,TL-JFZ:10,MH:12,HT-TL:15,GL-CG:17}.
  Interestingly, if we consider the linear SA algorithms over random signed networks with state-dependent noise,
from Theorems \ref{SAF_sc}, \ref{SAF_sn} and \ref{SAF_sn2} we can
obtain some results for the consensus or group consensus.

Assume the system contains $n$
agents. Each agent $i$ has a state $x_i(s)\in\real$ at time $s$
which can represent the opinion, social power or others, and is
updated according to the current state and the interaction from the
others.  In detail, for $1\leq i\leq n$ and $s\geq 0$, the state of
agent $i$ is updated by
\begin{multline}\label{model_consensus}
x_i(s+1)=(1-a(s))x_i(s)\\
+a(s)\sum_{j\in\mathcal{N}_i(s)} P_{ij}(s)\left[x_j(s)+f_{ji}(x(s))w_{ji}(s)\right],
\end{multline}
where $a(s)\geq 0$ is the gain at time $s$, $\mathcal{N}_i(s)$ is the neighbors of node $i$ at time $s$, $P_{ij}(s)$ is the weight of the edge
$(j,i)$ at time $s$, and $f_{ji}(x(s)) w_{ji}(s)$ is the noise of agent $i$ receiving information from agent $j$ at time $s$. Here
we consider the noise may be state-dependent which means that $f_{ji}(x(s))$ is a function of the state vector $x(s)$.
Let $P_{ij}(s)=0$ if $j\notin\mathcal{N}_i(s)$, and set
$$u_i(s):=\sum_{j\in\mathcal{N}_i(s)} P_{ij}(s)f_{ji}(x(s))w_{ji}(s),$$
then system (\ref{model_consensus}) can be rewritten as
 \begin{eqnarray*}
x(s+1)=(1-a(s))x(s)+a(s) \left[ P(s)x(s)+u(s)\right].
\end{eqnarray*}
If $P_{ij}(s)$ is a stationary stochastic process with uniformly bounded variance, and
$w_{ji}(s)$ is a zero-mean noise with uniformly bounded variance for any $x(s)$, $P_{ji}(s)$, and $j\in\mathcal{N}_i(s)$, then
(A1) is satisfied with $u={\mathbf{0}}_{n\times 1}$.

We say the subsets $S_1,\ldots,S_{r'} (r'\geq 1)$ is \emph{a partition
  of} $\{1,\ldots,n\}$ if $\emptyset\subset S_i \subseteq
\{1,\ldots,n\}$ for $1\leq i\leq r'$, $S_i \cap S_j=\emptyset$ for
$i\neq j$, and $\cup_{i=1}^{r'} S_i=\{1,\ldots,n\}$.  Following
\cite{JY-LW:10} with some modifications we introduce the definition
for group consensus:
\begin{definition}\label{def_consensus}
  Let the subsets $S_1,\ldots,S_{r'}$ be a partition of
  $\{1,\ldots,n\}$.  If $x(s)$ is mean-square convergent, and
  $\lim_{s\to\infty} \E|x_i(s)-x_j(s)|=0$ when $i$ and $j$ belong to a
  same subset, then we say $x(s)$ asymptotically reaches
  $\{S_i\}_{i=1}^{r'}$-group consensus in mean square.
\end{definition}

The group consensus turns to cluster consensus if different groups have different limit values \cite{YH-WL-TC:13}.

From Definition \ref{def_consensus} we can know that consensus is a
special case of the $\{S_i\}_{i=1}^{r'}$-group consensus with $r'=1$.
Before the statement of our results, we need to introduce some
notations and an assumption:

For a partition $S_1,\ldots,S_{r'}$ of $\{1,\ldots,n\}$, let $\mathds{1}^i\in\real^n (1\leq i\leq r')$ denote the column vector satisfying
$\mathds{1}_k^i=1$ if $k\in S_i$ and $\mathds{1}_k^i=0$ otherwise. A linear combination of $\{\mathds{1}^i\}_{i=1}^{r'}$ is $c_1 \mathds{1}^1+\ldots
+c_{r'} \mathds{1}^{r'}$ with $c_1,\ldots,c_{r'}\in\complex$ being constants.

\textbf{(A5)}
 Assume any eigenvalue of $P$ whose real part is $1$ equals $1$,
and the algebraic and geometric multiplicities of the eigenvalue $1$ equal $r\in[1,r']$,  and any right eigenvector of $P$ corresponding to the eigenvalue $1$
can be written as a linear combination of $\{\mathds{1}^i\}_{i=1}^{r'}$.

With Theorems \ref{SAF_sc}, \ref{SAF_sn} and Proposition \ref{SAF_sc} we obtain the following result:

\begin{theorem}\label{SAF_con_sn1}(Necessary and sufficient condition for group consensus with non-negative gains)
Consider the system (\ref{SAF_m1}) or (\ref{model_consensus}) satisfying (A1) with $u={\mathbf{0}}_{n\times 1}$ and (A2). Let $S_1,\ldots,S_{r'}$ be a partition of $\{1,\ldots,n\}$.
Then $x(s)$ asymptotically reaches $\{S_i\}_{i=1}^{r'}$-group consensus in mean square for any initial state
if and only if $\widetilde{\rho}_{\max}(P)<1$, or (A5) holds with $\widetilde{\rho}_{\max}(P)=1$.
\end{theorem}
\begin{IEEEproof}
Before proving our result, we introduce some notes first.
For any matrix $A\in\complex^{n\times n}$,  let $A_i$ and $A^{i}$ denote the $i$-th row and $i$-th column of $A$ respectively.
Set $A^{[i,j]}=(A^i,A^{i+1},\ldots,A^j)\in\complex^{n\times (j-i+1)}$.

We first consider the sufficient part. If $\widetilde{\rho}_{\max}(P)<1$, by Proposition \ref{SAF_sc} and the fact $u={\mathbf{0}}_{n\times 1}$ we obtain that $x(s)$ converges to ${\mathbf{0}}_{n\times 1}$ in mean square  for all initial states. Hence, the $\{S_i\}$-group consensus can be reached.

If (A5) holds with $\widetilde{\rho}_{\max}(P)=1$, which implies that (A3) holds together with the fact $u={\mathbf{0}}_{n\times 1}$.
Let $P=H^{-1} D H$, where $H$ is an invertible matrix, and $D$ is the Jordan normal form of $P$ with the same expression as (\ref{Jordan_temp}).
Then, by Theorem \ref{SAF_c1}, for any initial state there exist random variables $y_1,\ldots,y_r$ such that
in mean square
\begin{eqnarray}\label{con_sn1_1}
x(s) \rightarrow y_1 [H^{-1}]^1+\cdots+y_r [H^{-1}]^{r}~~\mbox{as}~~s\to\infty.
\end{eqnarray}
Also, from $P H^{-1}=D H^{-1}$ and (\ref{Jordan_temp}) we have
\begin{eqnarray}\label{con_sn1_2}
P [H^{-1}]^i = [H^{-1}]^i,~~1\leq i\leq r.
\end{eqnarray}
Hence, by  (\ref{con_sn1_1}) and (A5), there exist random variables $z_1,\ldots,z_{r'}$ such that in mean square
\begin{eqnarray*}\label{con_sn1_3}
x(s) \rightarrow z_1 \mathds{1}^1+\cdots+z_{r'} \mathds{1}^{r'}~~\mbox{as}~~s\to\infty,
\end{eqnarray*}
which implies that $x(s)$ asymptotically reaches $\{S_i\}_{i=1}^{r'}$-group consensus in mean square for any initial state.

Next we prove the necessary part.
Since $x(s)$ asymptotically reaches $\{S_i\}_{i=1}^{r'}$-group consensus in mean square for any initial state,
then, by Definition \ref{def_consensus}, $x(s)$ is mean-square convergent for any initial state. Hence, by Theorem \ref{SAF_sn},
we obtain that $\widetilde{\rho}_{\max}(P)<1$, or (A3) holds with $\widetilde{\rho}_{\max}(P)=1$.

It remains to show (A5) holds for the case when (A3) holds.
For any complex right eigenvector $\mathbf{a}+\mathbf{b}i\in\complex^n$ of $P$ corresponding to eigenvalue $1$, we have
$P \mathbf{a}=\mathbf{a}$ and  $P \mathbf{b}=\mathbf{b}$, which implies that $\mathbf{a}$ and $\mathbf{b}$ are real right eigenvectors
of $P$ corresponding to eigenvalue $1$. Thus, any complex right eigenvector of $P$ corresponding to the eigenvalue $1$ can be written as a linear combination of real right eigenvectors corresponding to the eigenvalue $1$. Also, from (\ref{Jordan_temp}) we have $PH^{-1}=H^{-1} D$ if and only if
(\ref{con_sn1_2}) and $P [H^{-1}]^{[r+1,n]}=[H^{-1}]^{[r+1,n]} \underline{D}$ hold.
Thus, we can choose suitable $H$ such that  $P=H^{-1} D H$ and $[H^{-1}]^1,\ldots,[H^{-1}]^{r}$ are real vectors.
 By Theorem \ref{SAF_c1}, we have
 \begin{eqnarray}\label{con_sn1_4}
\lim_{s\to\infty} \E x(s)=\sum_{i=1}^{r}H_i x(0) \cdot [H^{-1}]^i
\end{eqnarray}
Also, from $H H^{-1}=I_n$ we have $H_i [H^{-1}]^j$ equals $1$ if $i=j$ and $0$ otherwise.  If we choose
$x(0)=[H^{-1}]^i$ ($1\leq i\leq r$), by (\ref{con_sn1_4}) we have
$\lim_{s\to\infty} \E x(s)=[H^{-1}]^i$.
 Because for any initial state, $x(s)$ asymptotically reaches $\{S_i\}_{i=1}^{r'}$-group consensus in mean square, which implies
$Ex(s)$ also asymptotically reaches $\{S_i\}_{i=1}^{r'}$-group consensus, $[H^{-1}]^i (1\leq i\leq r)$ can be written as a linear combination of
$\{\mathds{1}^j\}_{j=1}^{r'}$. From the linear independence of $[H^{-1}]^1,\ldots,[H^{-1}]^{r}$
we have $r\leq r'$, and $[H^{-1}]^1,\ldots,[H^{-1}]^{r}$ is a basis of the eigenspace $\mathcal{R}^1$ which consists of all the right eigenvectors of $P$ corresponding to the eigenvalue $1$ and together with the zero vector. Hence, any vector in $\mathcal{R}^1$ can be written
as a linear combination of $[H^{-1}]^1,\ldots,[H^{-1}]^{r}$, and thus a linear combination of $\{\mathds{1}^j\}_{j=1}^{r'}$.
\end{IEEEproof}

Similar to Theorem \ref{SAF_con_sn1} we have the following theorem:
\begin{theorem}\label{SAF_con_sn2}(Necessary and sufficient condition for group consensus with non-positive gains)
Consider the system (\ref{SAF_m1})  or (\ref{model_consensus}) satisfying (A1) with $u={\mathbf{0}}_{n\times 1}$ and (A2').
Let $S_1,\ldots,S_{r'}$ be a partition of $\{1,\ldots,n\}$.
Then $x(s)$ asymptotically reaches $\{S_i\}_{i=1}^{r'}$-group consensus in mean square for any initial state
if and only if $\widetilde{\rho}_{\min}(P)>1$, or (A5) holds with $\widetilde{\rho}_{\min}(P)=1$.
\end{theorem}

By Theorems \ref{SAF_con_sn1} and \ref{SAF_con_sn2} with $r'=1$,  we immediately obtain the following two corollaries for consensus:

\begin{corollary}\label{SAF_con_sn3}
Consider the system (\ref{SAF_m1}) or (\ref{model_consensus}) satisfying (A1) with $u={\mathbf{0}}_{n\times 1}$ and (A2).
Then $x(s)$ asymptotically reaches consensus in mean square for any initial state
if and only if one of the following condition holds:\\
i) $\widetilde{\rho}_{\max}(P)<1$;\\
ii) The sum of  of each row of $P$ equals $1$, and $P$
has $n-1$ eigenvalues whose real parts are all less than $1$.
\end{corollary}

\begin{corollary}\label{SAF_con_sn4}
Consider the system (\ref{SAF_m1}) or (\ref{model_consensus}) satisfying (A1) with $u={\mathbf{0}}_{n\times 1}$ and (A2').
Then $x(s)$ asymptotically reaches consensus in mean square for any initial state
if and only if one of the following condition holds:\\
i) $\widetilde{\rho}_{\min}(P)>1$;\\
ii) The sum of  of each row of $P$ equals $1$, and $P$
has $n-1$ eigenvalues whose real parts are all bigger than $1$.
\end{corollary}

The communication topology is an important aspect in the research of multi-agent systems consensus.
In fact, our result can also give some topology conditions of consensus for some special $P$.
We first introduce some definitions concerning graphs.
For a matrix $A\in \real^{n\times n}$ with  $A_{ij}\geq 0$ for $j\neq i$.
let $\mathcal{V}=\{1,2,\ldots n\}$ denote the set of nodes, and $\mathcal{E}$ denote the set of edges
where an ordered pair  $(j,i)\in \mathcal{E}$  if and only if $A_{ij}>0$. The digraph associated with $A$ is defined by
$\mathcal{G}=\{\mathcal{V},\mathcal{E}\}$.
A sequence $(i_1, i_2), (i_2, i_3), \ldots, (i_{k-1}, i_k)$  of edges is called
a directed path from node $i_1$ to node $i_k$. $\mathcal{G}$ contains a directed spanning tree if
there exists a root node $i$ such that $i$ has a directed path to $j$  for any node $j\neq i$.

We need the following lemma in our results.
\begin{lemma}[Lemma 3.3 in \cite{WR-RWB:05}]\label{Ren05}
Given a matrix $A\in \real^{n\times n}$, where for any $i\in\mathcal{V}$, $A_{ii}\leq 0, A_{ij}\geq 0$ for $j\neq i$, and $\sum_{j=1}^n A_{ij}=0$,  then $A$ has at least
one zero eigenvalue and all of the non-zero eigenvalues have negative real parts. Furthermore, $A$ has exactly one zero eigenvalue
if and only if the directed graph associated with $A$ contains a directed spanning
tree.
\end{lemma}

From Corollary \ref{SAF_con_sn3} and Lemma \ref{Ren05} we have the following result.

\begin{corollary}\label{SAF_con_sn5}
Consider the system (\ref{SAF_m1}) or (\ref{model_consensus}) satisfying (A1) and (A2).
Assume that $P$ is a row-stochastic matrix and $u={\mathbf{0}}_{n\times 1}$. Then $x(s)$ asymptotically reaches consensus in mean square for any initial state
if and only if the digraph associated with $P$ contains a directed spanning
tree.
\end{corollary}
\begin{IEEEproof}
Let $A=P-I_n$ and `$\leftrightarrow$' denote the `if and only if'.  The digraph associated with $P$ contains a directed spanning
tree $\leftrightarrow$  the digraph associated with $A$ contains a directed spanning
tree $\stackrel{\underleftrightarrow{Lemma ~\ref{Ren05}}}{ }$
$A$ has exactly one zero eigenvalue, and all the non-zero
 eigenvalues have negative real parts $\leftrightarrow$
 $P$
has $n-1$ eigenvalues whose real parts are all less than $1$ $\stackrel{\underleftrightarrow{Corollary~\ref{SAF_con_sn3}}}{ }$
 $x(s)$ asymptotically reaches consensus in mean square for any initial state, where the last two `$\leftrightarrow$'
uses the hypothesis that $P$ is a row-stochastic matrix which has at least one eigenvalue that is equal to $1$.
\end{IEEEproof}

Corollary \ref{SAF_con_sn5} coincides with the consensus condition for the continuous-time consensus protocol with time-invariant interaction topology
(Theorem 3.8 in \cite{WR-RWB:05}).

If $P$ is not a row-stochastic matrix, the consensus may be also reached.
For example, let

\begin{eqnarray}\label{eq:consensusmatrix}
P:=
\begin{bmatrix}
    0.5 &0.3 &0 &0.3& -0.1   \\
    -0.1& 0.3 &0.3 &0 &0.5   \\
     0& 0.2 &0.4& 0.5& -0.1  \\
     0.1& 0& 0.6 &0.4 &-0.1 \\
     0.1& -0.1& 0.1& 0.3& 0.6
\end{bmatrix}.
\end{eqnarray}
The eigenvalues of $P$ are $1,  0.5708, -0.2346,0.4319+0.3270i,0.4319-0.3270i$. By Corollary~\ref{SAF_con_sn3} $x(s)$  asymptotically reaches consensus in mean square.

Different from consensus, the group consensus does not require that the sum of each row of $P$ equals $1$. For example,
if
\begin{eqnarray}\label{eq:gconsensusmatrix}
P=
\begin{bmatrix}
    0.3 & 0.5 & 0.5 & -0.4 \\
    0.5 & 0.3 & -0.4 & 0.5 \\
     -0.1 & 0.5 & 0.4 & 0.4\\
     0.5 & -0.1 & 0.4 & 0.4
\end{bmatrix},
\end{eqnarray}
then $P[1,1,2,2]^\top=[1,1,2,2]^\top$, and the eigenvalues of $P$ are $1,0.6,-0.1+0.728i,-0.1-0.728i$.
Let $S_1=\{1,2\}$ and $S_2=\{3,4\}$, by Theorem \ref{SAF_con_sn1}
$x(s)$ can asymptotically reach $\{S_1,S_2\}$-group consensus in mean square for any initial state.

\begin{figure}
\begin{center}
\subfigure[Consensus]{
\includegraphics[scale=0.4]{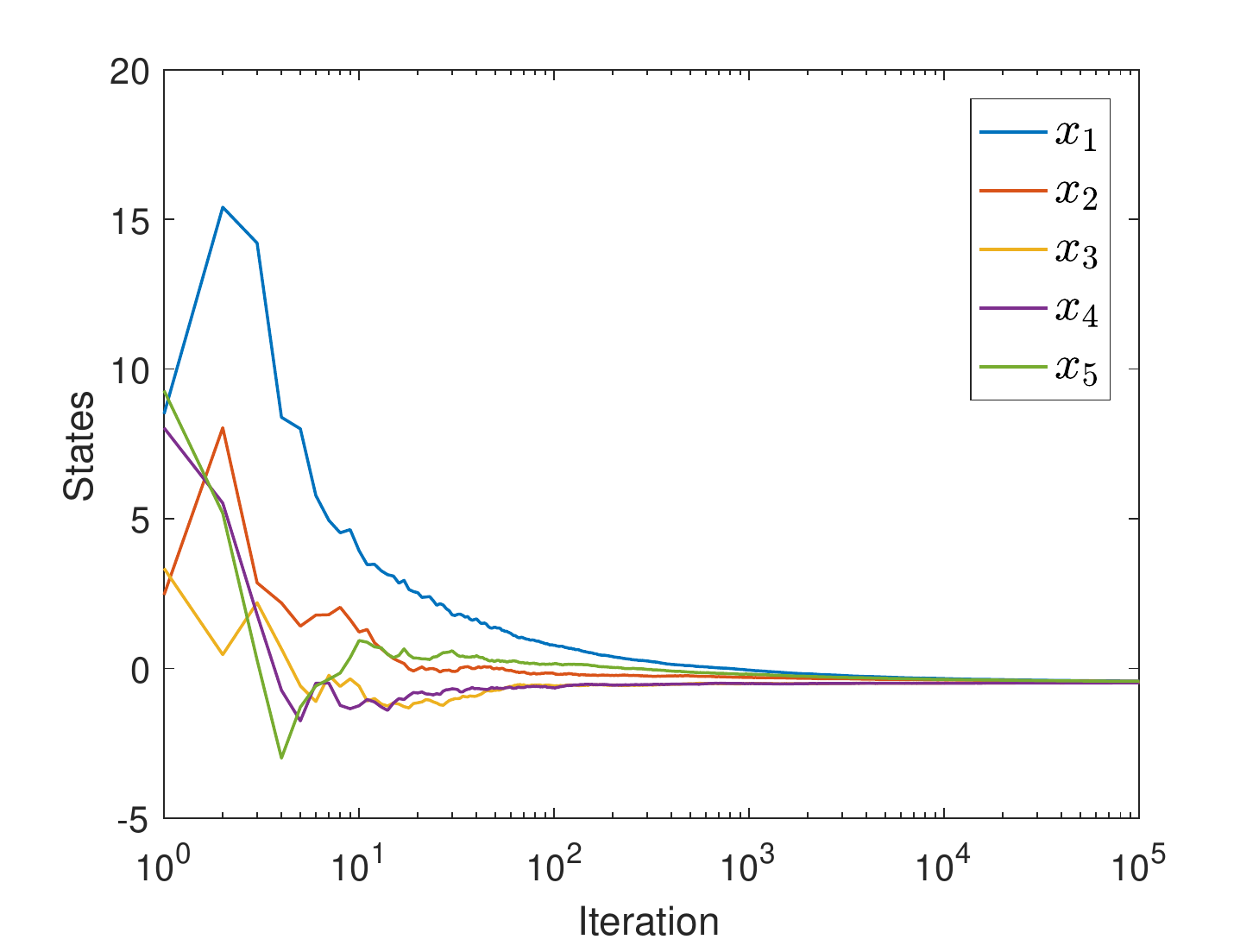}}
\subfigure[ Group consensus]{
\includegraphics[scale=0.4]{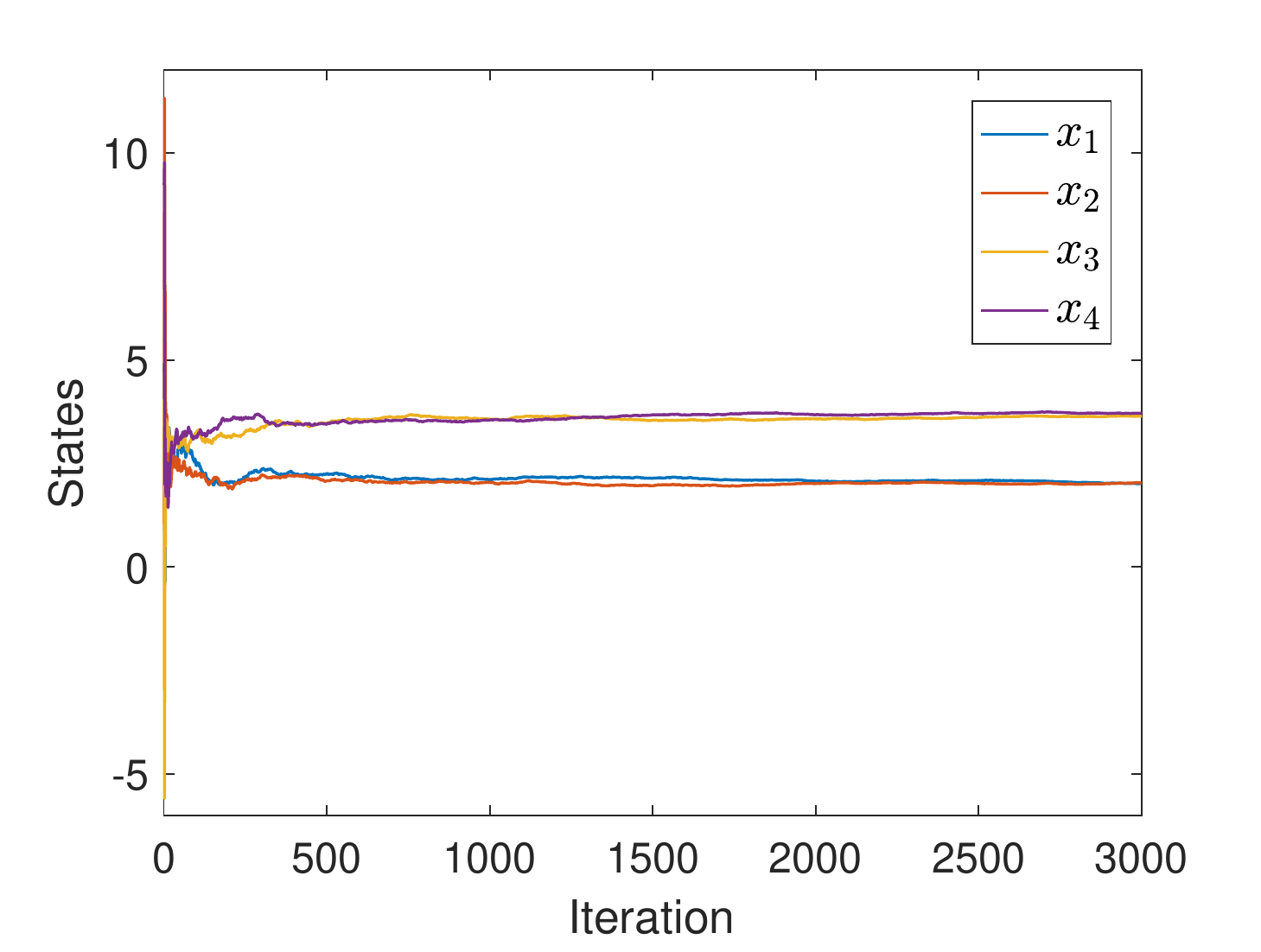}}
\caption{Consensus and group consensus under linear SA protocol (\ref{SAF_m1})}\label{fig:consensusunderSA}
\end{center}
\end{figure}

In the following, we simulate system (\ref{SAF_m1}) to show consensus and group consensus using $P$ matrices in \eqref{eq:consensusmatrix} and \eqref{eq:gconsensusmatrix} respectively. For $s\geq 0$, $P(s)$ and $u(s)$ are generated by i.i.d. matrix and vector
with mean $P$ and ${\mathbf{0}}_{n\times 1}$ respectively.
 We set the gain function $a(s)=\frac{1}{s}$. From Fig. \ref{fig:consensusunderSA}, we can see that consensus and group consensus are reached as guaranteed by Corollary \ref{SAF_con_sn3} and Theorem \ref{SAF_con_sn1}, respectively.

\subsection{An extension to multidimensional linear SA algorithms}

Our results in Section \ref{Main_results} can be extended to multidimensional linear SA algorithms in which
the state of each agent is a $m$-dimensional vector.
The dynamics is, for all $s\geq 0$
\begin{equation}\label{MSAF}
X(s+1)=(1-a(s))X(s)+a(s)\big[P(s)X(s)C^\top(s)+U(s)\big],
\end{equation}
where $X(s)\in\real^{n\times m}$ is the state matrix, $P(s)\in\real^{n\times n}$ is still an interaction matrix,
$C\in\real^{m\times m}$ is an interdependency matrix,
 and $U(s)\in\real^{n\times m}$ is
an input matrix.

The system (\ref{MSAF}) can be transformed to one dimensional system (\ref{SAF_m1}) by the following way:


Given a pair of matrices $A\in\real^{n\times m}$, $B\in\real^{p\times q}$,
 their Kronecker product is defined by
 \begin{eqnarray*}
A\otimes B=
\begin{bmatrix}
    A_{11}B & \cdots  & A_{1m} B    \\
     \cdots  & \cdots & \cdots   \\
     A_{n1}B  &  \cdots  & A_{nm} B
\end{bmatrix}\in\real^{np\times mq}.
\end{eqnarray*}
Let $Q(s):=P(s)\otimes C(s)$.
From (\ref{MSAF}) we have
\begin{eqnarray}\label{MSAF_sc_1}
&&X_{ij}(s+1)=(1-a(s))X_{ij}(s)\\
&&~~+a(s)\Big[\sum_{k_1,k_2} P_{ik_1}(s)X_{k_1,k_2}(s)C_{jk_2}(s)+U_{ij}(s)\Big].\nonumber
\end{eqnarray}
for any $s\geq 0$, $1\leq i\leq n$, and $1\leq j\leq m$.
Let
\begin{eqnarray*}\label{MSAF_sc_2}
y(s):=(X_{11}(s),\ldots,X_{1m}(s),\ldots,X_{n1}(s),\ldots,X_{nm}(s))^\top
\end{eqnarray*}
and
\begin{eqnarray*}\label{MSAF_sc_3}
v(s):=~(U_{11}(s),\ldots,U_{1m}(s),\ldots,U_{n1}(s),\ldots,U_{nm}(s))^\top
\end{eqnarray*}
be the vector in $\real^{nm}$ transformed from the matrices $X(s)$ and $U(s)$ respectively.
By (\ref{MSAF_sc_1}) we have

\begin{align*}
y&_{(i-1)m+j}(s+1)=X_{ij}(s+1)\\
&=(1-a(s))X_{ij}(s)\\
&\quad+a(s)\Big[\sum_{k_1,k_2} P_{ik_1}(s)X_{k_1,k_2}(s)C_{jk_2}(s)+U_{ij}(s)\Big]\\
&=(1-a(s))y_{(i-1)m+j}(s)+\\
&\quad a(s)\Big[\sum_{k_1,k_2}Q_{(i-1)m+j,(k_1-1)m+k_2}(s)y_{(k_1-1)m+k_2}(s)\\
&\quad\quad\quad+v_{(i-1)m+j}(s)\Big],
\end{align*}
which implies
\begin{eqnarray*}\label{MSAF_sc_5}
y(s+1)=(1-a(s))y(s)+a(s)[Q(s)y(s)+v(s)].
\end{eqnarray*}

The system (\ref{MSAF_sc_5}) has the same form as the system (\ref{SAF_m1}), so the results in
Section \ref{Main_results} can be applied to the multidimensional linear SA algorithms.

\subsection{SA Friedkin-Johnsen model over time-varying interaction network}

The Friedkin-Johnsen (FJ) model proposed by \cite{NEF-ECJ:99}
considers a community of $n$ social actors (or agents) whose opinion
column vector is $x(s)=(x_1(s),\dots, x_n(s))^\top \in\real^{n}$ at time
$s$. The FJ model also contains a row-stochastic matrix of
interpersonal influences $P\in\real^{n\times n}$ and a diagonal matrix
of actors' susceptibilities to the social influence
$\Lambda\in\real^{n\times n}$ with ${\mathbf{0}}_{n\times n}\leq \Lambda \leq I_n$.  The state
of the FJ model is updated by
\begin{eqnarray}\label{FJ}
x(s+1)=\Lambda P x(s)+(I_n-\Lambda)x(0),~~s=0,1,\ldots.
\end{eqnarray}
By \cite{SEP-AVP-RT-NEF:17}, if ${\mathbf{0}}_{n\times n}\leq \Lambda<I_n$, then
\begin{eqnarray}\label{FJ_r1}
\lim_{s\to\infty}x(s)=(I_n-\Lambda P)^{-1}(I_n-\Lambda)x(0).
\end{eqnarray}
However, if the interpersonal influences are affected by noise, then the system (\ref{FJ}) may not converge.


The FJ model (\ref{FJ}) was extended to the multidimensional case in
\cite{SEP-AVP-RT-NEF:17,NEF-AVP-RT-SEP:16}.  The multidimensional FJ model still
contains $n$ individuals, but each individuals has beliefs on $m$ truth
statements.  Let $X(s)\in\real^{n\times m}$ be the matrix of $n$
individuals' beliefs on $m$ truth statements at time $s$.  Following
\cite{NEF-AVP-RT-SEP:16}, it is updated by
\begin{eqnarray}\label{MFJ}
X(s+1)=\Lambda P X(s) C^\top +(I_n-\Lambda)X(0)
\end{eqnarray}
for $s=0,1,\ldots,$ where $\Lambda,P\in\real^{n\times n}$ are the same
matrices in (\ref{FJ}), and $C\in\real^{m\times m}$ is a
row-stochastic matrix of interdependencies among the $m$ truth
statements. The convergence of system (\ref{MFJ}) has been analyzed in
\cite{SEP-AVP-RT-NEF:17}. Similar to (\ref{FJ}) it is easy to see that
if system (\ref{MFJ}) is affected by noise, then it will not
converge.  We will adopt the stochastic-approximation method to smooth
the effects of the noise.
\begin{proposition}\label{SFJ_r2}
Consider the system
\begin{multline}\label{SFJ_r2_00}
X(s+1)=(1-a(s))X(s)+a(s)[\Lambda(s) P(s) X(s)C(s)^\top \\
+(I_n-\Lambda(s))X(0)],
\end{multline}
for $s=0,1,\ldots$, where $\Lambda(s)\in\real^{n\times n}$, $P(s)\in \real^{n\times n}$ and $C(s)\in \real^{m\times m}$ are independent matrix sequence
with invariant expectation $\Lambda$, $P$, and $C$ respectively. Assume $E\|\Lambda(s)\|_2^2$,  $E\|P(s)\|_2^2$, and  $E\|C(s)\|_2^2$ are uniformly bounded.
 Suppose
$P$ and $C$ are row-stochastic matrix, and ${\mathbf{0}}_{n\times n}\leq \Lambda<I_n$, and
 the gain function $a(s)$ satisfies (A2).
Then for any initial state, $X(s)$ converges to $X^*$ in mean square,
where $X^*$ is the unique solution of the equation
\begin{eqnarray}\label{SFJ_r2_01}
X=\Lambda P  X C^\top+(I_n-\Lambda)X(0).
\end{eqnarray}
\end{proposition}
\begin{IEEEproof}
Since $P$ and $C$ are row-stochastic matrices, $P\otimes C$ is still a row-stochastic matrix. Together with the condition that ${\mathbf{0}}_{n\times n}\leq \Lambda <I_n$,
we have that the sum of each row of $(\Lambda P)\otimes C$ is less than $1$. Thus, using the Ger\v{s}gorin Disk Theorem we obtain
 $\widetilde{\rho}_{\max}((\Lambda P)\otimes C)<1$.
Let $Q:=(\Lambda P)\otimes C$, $U(s):=(I_n-\Lambda(s))X(0)$,
\begin{eqnarray*}
y(s):=(X_{11}(s),\ldots,X_{1m}(s),\ldots,X_{n1}(s),\ldots,X_{nm}(s))^\top,
\end{eqnarray*}
\begin{eqnarray*}
v(s):=(U_{11}(s),\ldots,U_{1m}(s),\ldots,U_{n1}(s),\ldots,U_{nm}(s))^\top,
\end{eqnarray*}
and $v:=Ev(s)$.
 By Proposition \ref{SAF_sc} and the transformation from (\ref{MSAF}) to (\ref{MSAF_sc_5}), we obtain that
$y(s)$ converges to $(I_{mn}-Q)^{-1}v$ in mean square.

It remains to discuss the relation between $(I_{mn}-Q)^{-1} v$ and $X^*$.
Let
\begin{eqnarray*}\label{MSAF_sc_6}
y^*:=(X_{11}^*,\ldots,X_{1m}^*,\ldots,X_{n1}^*,\ldots,X_{nm}^*)^\top\in\real^{nm}.
\end{eqnarray*}
By (\ref{SFJ_r2_01}),  similar to (\ref{MSAF_sc_5}) we have $y^*=Qy^*+v$, which has a unique solution
$y^*=(I_{mn}-Q)^{-1} v$ since $I_{mn}-Q$ is an invertible matrix by $\widetilde{\rho}_{\max}(Q)<1$.
Thus, with the fact that $y(s)$ converges to $(I_{mn}-Q)^{-1} v$ in mean square we obtain that $X(s)$ converges to $X^*$ in mean square.
\end{IEEEproof}

\begin{remark}
According to Theorem \ref{SAF_c1} and Proposition \ref{SAF_sc}, the conditions of $\Lambda(s),P(s)$ and $C(s)$ in Proposition \ref{SFJ_r2} can be further relaxed for convergence, such as $P$ and $C$ are not row-stochastic matrices, and ${\mathbf{0}}_{n\times n}\leq \Lambda<I_n$ may be extended to $\Lambda<{\mathbf{0}}_{n\times n}$ or $\Lambda\geq I_n$.
\end{remark}

\renewcommand{\thesection}{\Roman{section}}
\section{Conclusion}\label{sec:conclusion}
\renewcommand{\thesection}{\arabic{section}}
In this paper, we study a time-varying linear dynamical system,
where the state of the system features persistent oscillation and does
not converge. We consider a stochastic approximation-based approach and obtain necessary and sufficient conditions to guarantee
mean-square convergence. Our theoretical results largely extend the
conditions on the spectrum of the expectation of the system matrix and
thus can be applied in a much broader range of applications. We also
derived the convergence rate of the system. To illustrate the
theoretical results, we applied them in two different applications:
group consensus in multi-agent systems and FJ model with
time-varying interactions in social networks.

This work leaves various problems for future research. First, the
system matrix and input are assumed to have constant expectations in
this paper. However, it would be more interesting, yet challenging, to
study systems with time-varying expectation of the system matrix and
input. Second, we only considered  linear dynamical systems in this
paper. How and whether the proposed framework can be extended to
non-linear system are important and intriguing questions. Finally, we
have illustrated our results in two different application scenarios;
there are other possible applications such as gossip algorithms for
consensus.

\appendices

\section{}\label{App_lemmas}

\begin{lemma}\label{lguo}
Suppose the non-negative real number sequence $\{y_s\}_{s\geq 1}$ satisfies
\begin{equation}\label{lguo_00}
y_{s+1}\leq (1-a_s)y_s+b_s,
\end{equation}
where  $b_s\geq 0$ and $a_s\in[0,1)$ are real numbers.
If $\sum_{s=1}^{\infty}a_s =\infty$ and $\lim_{s\to\infty} b_s/a_s=0$, then
$\lim_{s\to\infty}y_s=0$ for  any $y_1\geq 0$.
\end{lemma}
\begin{IEEEproof}
Repeating (\ref{lguo_00}) we obtain
\begin{eqnarray*}\label{lguo_1}
\begin{aligned}
y_{s+1}&\leq y(1) \prod_{t=1}^{s}(1-a_t)+\sum_{i=1}^s b_i \prod_{t=i+1}^s (1-a_t).
\end{aligned}
\end{eqnarray*}
Here we define $\prod_{t=i}^s (\cdot):=1$ when $i>s$. From the hypothesis $\sum_{t=1}^{\infty} a_t=\infty$
we have $\prod_{t=1}^{\infty}(1-a_t)=0$. Thus,
to obtain $\lim_{s\to\infty}y_s=0$ we just need to prove that
\begin{eqnarray}\label{lguo_2}
\lim_{s\to\infty}\sum_{i=1}^s b_i \prod_{t=i+1}^s (1-a_t)=0.
\end{eqnarray}
Since $\lim_{s\to\infty} b_s/a_s=0$, for any real number $\varepsilon>0$, there exists an integer $s^*>0$ such that
$b_s\leq \varepsilon a_s$  when $s\geq s^*$. Thus,
\begin{align}\label{lguo_3}
\sum_{i=1}^s &b_i \prod_{t=i+1}^s (1-a_t)\\
&\leq \sum_{i=1}^{s^*-1} b_i \prod_{t=i+1}^s (1-a_t)+\sum_{i=s^*}^{s} \varepsilon a_i \prod_{t=i+1}^s (1-a_t)\nonumber\\
&=\sum_{i=1}^{s^*-1} b_i \prod_{t=i+1}^s (1-a_t)+ \varepsilon\bigg( 1-\prod_{t=s^*}^s (1-a_t)\bigg)\nonumber\\
&\rightarrow \varepsilon ~~\mbox{as}~~s\to\infty,\nonumber
\end{align}
 where the first equality uses the classic equality
 \begin{eqnarray}\label{lguo_4}
\sum_{t=s^*}^{s} c_t \prod_{k=t+1}^{s} (1-c_k)=1-\prod_{t=s^*}^{s} (1-c_t)
\end{eqnarray}
with $\{c_t\}$ being any complex numbers, which can be obtained by induction. Here we define $\prod_{k=s_1}^{s_2}(\cdot)=1$ if $s_2<s_1$.
 Let $\varepsilon$ decrease to $0$, then (\ref{lguo_3}) is followed by (\ref{lguo_2}).
\end{IEEEproof}

\section{Proof of Theorem \ref{SAF_r1}}\label{proof_SAF_r1}
We prove this theorem under the following three cases:

\textbf{Case I:} $\widetilde{\rho}_{\max}(P)<1$.  Define $\theta(s)$,
$A$ and $A_2$ as in the proof of Theorem \ref{SAF_c1} but with $r=0$.
Set $V(\theta):=\theta^* A \theta$ for any $\theta\in\complex^n$,
where $\theta^*$ denotes the conjugate transpose of $\theta$.  We
remark that $A_2=A\in\complex^{n\times n}$ under the case $r=0$, so
that, by (\ref{SAF_c1_16}), we have
\begin{align}\label{SAF_cr_1}
\E[&V(\theta(s+1))]\\
&\leq \Big(1-\frac{\alpha}{\rho(A)(s+\beta)^{\gamma}}\Big) \Big) \E[V(\theta(s))]+O\Big(\frac{1}{(s+\beta)^{2\gamma}} \Big).\nonumber
\end{align}
Set
 \begin{eqnarray*}\label{SAF_cr_2}
\Phi(s,i):=\prod_{k=i}^s  \Big(1-\frac{\alpha}{\rho(A)(k+\beta)^{\gamma}}\Big)
\end{eqnarray*}
 and define $\prod_{k=i}^s (\cdot):=1$ if $s<i$. We  compute
 \begin{align}\label{SAF_cr_3}
\Phi(s,i)& =O \bigg( \exp\Big[ \sum_{k=i}^s  -\frac{\alpha}{\rho(A)(k+\beta)^{\gamma}} \Big]\bigg)\\
&=O \bigg(\exp\Big( \int_{i}^s -\frac{\alpha}{\rho(A)(k+\beta)^{\gamma}}dk \Big) \bigg)\nonumber\\
&=
\left\{
\begin{array}{ll}
O \Big(\big(\frac{s+\beta}{i+\beta}\big)^{-\alpha/\rho(A)}\Big),  \mbox{ if } \gamma=1, \nonumber\\
O\Big(\exp\big(\frac{-\alpha}{(1-\gamma)\rho(A)}[(s+\beta)^{1-\gamma}-(i+\beta)^{1-\gamma}]\big)\Big),\nonumber\\
~~~~~~~~~~~~~~~~~~~~~~~~\mbox{ if } \frac{1}{2}<\gamma<1.\nonumber
\end{array}
\right.
\end{align}
Also, using (\ref{SAF_cr_1}) repeatedly we obtain
\begin{align}\label{SAF_cr_4}
\E[&V(\theta(s+1))]\\
&\leq \Phi(s,0)\E[V(\theta(0))]+\sum_{i=0}^s \Phi(s,i+1)O\Big( \frac{1}{(i+\beta)^{2\gamma}}\Big).\nonumber
\end{align}
Assume $\alpha\geq \rho(A)$. We first consider the case that
$\gamma=1$. From (\ref{SAF_cr_3}) and (\ref{SAF_cr_4}) we have
\begin{equation}\label{SAF_cr_5}
\E[V(\theta(s+1))]
= o\Big(\frac{1}{s}\Big)+O\bigg(\sum_{i=0}^s \frac{(s+\beta)^{-\alpha/\rho(A)}}{(i+\beta)^{2-\frac{\alpha}{\rho(A)}}}\bigg)=O\Big(\frac{1}{s}\Big).
\end{equation}
For the case when $\gamma\in(\frac{1}{2},1)$, we take
$b=\frac{\alpha}{(1-\gamma)\rho(A)}$, and from (\ref{SAF_cr_3}) and
(\ref{SAF_cr_4}) we can obtain
\begin{align}\label{SAF_cr_6}
\E[&V(\theta(s+1))]\nonumber\\
&=e^{-b(s+\beta)^{1-\gamma}}\cdot O\bigg(1+ \sum_{i=0}^s \frac{e^{b(i+\beta)^{1-\gamma}}}{(i+\beta)^{2\gamma}}\bigg)\nonumber\\
&=e^{-b(s+\beta)^{1-\gamma}}\cdot O\bigg(\sum_{i=0}^s \sum_{k=0}^{\infty}  \frac{b^k (i+\beta)^{(1-\gamma)k-2\gamma}}{k!}\bigg)\nonumber\\
&=e^{-b(s+\beta)^{1-\gamma}}\cdot O\bigg(\sum_{k=0}^{\infty} \frac{b^k}{k!} \sum_{i=0}^s  (i+\beta)^{(1-\gamma)k-2\gamma}\bigg)\nonumber\\
&=e^{-b(s+\beta)^{1-\gamma}}\cdot O\bigg(\sum_{k=0}^{\infty} \frac{b^k (s+\beta)^{(1-\gamma)k-2\gamma+1}}{k![(1-\gamma)k-2\gamma+1]} \bigg)\nonumber\\
&=\frac{e^{-b(s+\beta)^{1-\gamma}}}{(s+\beta)^{\gamma}}\cdot O\bigg(\sum_{k=0}^{\infty} \frac{b^{k+1} (s+\beta)^{(1-\gamma)(k+1)}}{(k+1)!} \bigg)\nonumber\\
&=O\big(s^{-\gamma} \big).
\end{align}
By (\ref{SAF_cr_5}) and (\ref{SAF_cr_6}), we have
$\E[V(\theta(s))]=O(s^{-\gamma})$ for $\frac{1}{2}<\gamma\leq
1$. Combing this with the definition of $\theta(s)$ yields our result.

\textbf{Case II:} $\widetilde{\rho}_{\max}(P)=1$.  Let $\theta(s)$,
$\underline{\theta}(s)$, $\bar{\theta}(s)$, $\bar{\theta}(\infty)$,
$H$, $y$ and $z$ be the same variables as in the proof of Theorem
\ref{SAF_c1}.  With (\ref{SAF_c1_16}) and following the similar
process from (\ref{SAF_cr_1}) to (\ref{SAF_cr_6}), we have
$\E\|\underline{\theta}(s)\|_2^2=O(s^{-\gamma})$. Also, from
(\ref{SAF_c1_20}) we have
\begin{align*}
  \E\|\bar{\theta}(\infty)-\bar{\theta}(s)\|_2^2&=O\Big(\sum_{k=s}^{\infty}a^2(k)\Big)=O\Big(\sum_{k=s}^{\infty}a^2(k)\Big) \\
  &=O\Big(\sum_{k=s}^{\infty}\frac{1}{(s+\beta)^{2\gamma}}\Big)=O\Big(\frac{1}{s^{2\gamma-1}}\Big).
\end{align*}
Since $x(s)=H^{-1}[\theta(s)+z]$ and $H^{-1}y$ is a mean square limit
of $x(s)$, the arguments above imply
\begin{eqnarray*}\label{SAF_cr_8}
  \begin{aligned}
    \E\|x(s)-H^{-1}y\|_2^2&=\max\big\{O\big(s^{-\gamma}\big), O\big(s^{1-2\gamma}\big) \big\}\\
    &= O\big(s^{1-2\gamma}\big).
  \end{aligned}
\end{eqnarray*}

\textbf{Case III:} $\widetilde{\rho}_{\min}(P)\geq 1$. The protocol
(\ref{SAF_m1}) is written as
 \begin{eqnarray*}\label{SAF_cr_9}
x(s+1)=x(s)+\frac{\alpha}{(s+\beta)^{\gamma}}[(I_n-P(s)) x(s)-u(s)].
\end{eqnarray*}
Because $\widetilde{\rho}_{\max}(I_n-P)\leq 0$, arguments similar to
that for Cases I) and II) yield our result.

\section{Proof of Theorem \ref{SAF_nca}}\label{proof_nca}

i) As same as Subsection \ref{sf_2}, the Jordan normal
form of $H$ is
$$D=\mbox{diag}(J_1,\ldots,J_k)=H P H^{-1}.$$ We also set $y(s):=H
x(s)$, $v(s):=H u(s)$, $D(s):=H P(s) H^{-1}$, $D=\E D(s)=HPH^{-1}$, and
$v=Ev(s)=Hu$.  By (\ref{SAF_c1_5}) and (A1) we have
\begin{align}\label{SAF_nn_1}
\E y(s+1)&=\E\big[\E[y(s+1)\,|\,y(s)]\big]\nonumber\\
&=\E y(s)+a(s)[(D-I_n)Ey(s)+v].
\end{align}
Let $B(s):=I_n+a(s)(D-I_n)$. Using (\ref{SAF_nn_1}) repeatedly we obtain
\begin{equation}\label{SAF_nn_2}
\E[y(s+1)]
= B(s)\cdots B(0)y(0) +\sum_{t=0}^s a(t) B(s)\cdots B(t+1)v.
\end{equation}
We will continue the proof under the following two
cases:\\ \textbf{Case I}: $\widetilde{\rho}_{\max}(P)>1$.  Without
loss of generality we assume $\mbox{Re}(\lambda_{1}(P))>1$.  Let $J_1$
be a Jordan block in $D$ corresponding to $\lambda_{1}(P)$.  Let $m_1$
be the row index of $D$ corresponding to the last line of $J_1$, i.e.,
\begin{eqnarray}\label{SAF_nn_3}
D_{m_1}=(0,\ldots,0,\lambda_{1}(P),0,\ldots,0).
\end{eqnarray}
Then by (\ref{SAF_nn_2})
\begin{align}\label{SAF_nn_4}
\E[&y_{m_1}(s+1)]\nonumber\\
&= y_{m_1}(0)\prod_{t=0}^s [1-a(i)[1-\lambda_{1}(P)]]+\frac{v_{m_1}}{1-\lambda_{1}(P)}\nonumber\\
&\quad \times \sum_{t=0}^s a(t)[1-\lambda_{1}(P)] \prod_{k=t+1}^s (1-a(k)[1-\lambda_{1}(P)])\nonumber\\
&= y_{m_1}(0)\prod_{t=0}^s \left(1-a(t)[1-\lambda_{1}(P)]\right)+\frac{v_{m_1}}{1-\lambda_{1}(P)}\nonumber\\
&\quad \times \Big(1- \prod_{t=0}^s \left(1-a(t)[1-\lambda_{1}(P)]\right)\Big),
\end{align}
where the last equality uses the equality (\ref{lguo_4}).
Since $\sum_{s} a(s)=\infty$,
\begin{multline}\nonumber
\prod_{t=0}^{\infty} |1-a(t)[1-\lambda_{1}(P)]|^2\\
\geq \prod_{t=0}^{\infty} \big\{1+2a(t)[\mbox{Re}(\lambda_{1}(P))-1]\}=\infty.
\end{multline}
Hence, from (\ref{SAF_nn_4}), if $y_{m_1}(0)\neq \frac{v_{m_1}}{1-\lambda_{1}(P)}$, then
\begin{eqnarray}\label{SAF_nn_5_1}
\lim_{s\to\infty} |\E[y_{m_1}(s)]|=\infty,
\end{eqnarray}
 which implies $\lim_{s\to\infty}\E\|x(s)\|_2^2=\infty$.\\
\textbf{Case II}:  $\widetilde{\rho}_{\max}(P)=1$. Under this case we consider the following three
situations:\\
(a) There is an eigenvalue $\lambda_j(P)=1+\mbox{Im}(\lambda_{j}(P))i$ with $\mbox{Im}(\lambda_{j}(P))\neq 0$, where
$\mbox{Im}(\lambda_{j}(P))$ denotes the imaginary part of $\lambda_{j}(P)$.
Similar to (\ref{SAF_nn_3}), we can choose a row $D_{j'}$ of $D$ which is equal to
$(0,\ldots,0,\lambda_{j}(P),0,\ldots,0)$. Similar to (\ref{SAF_nn_4}), we have
\begin{align}
\E[y_{j'}(s+1)]&= y_{j'}(0)\prod_{t=0}^s [1-a(t)[1-\lambda_{j}(P)]]\nonumber\\
& +\frac{v_{j'}}{1-\lambda_{j}(P)}\cdot \Big(1- \prod_{t=0}^s [1-a(t)[1-\lambda_{j}(P)]]\Big)\label{SAF_nn_6}.
\end{align}
 We write
\begin{align*}
1-a(t)[1-\lambda_{j}(P)]&=1+a(t)\mbox{Im}(\lambda_{j}(P))i \\
&=r_t e^{i \varphi_t}=r_t(\cos \varphi_t+i\sin\varphi_t),
\end{align*}
where $r_t=\sqrt{1+a^2(t)\mbox{Im}^2(\lambda_{j}(P))}$ and
\begin{eqnarray}\label{SAF_nn_8}
\begin{aligned}
\varphi_t&=\arctan [a(t)\mbox{Im}(\lambda_{j}(P))]\\
&=a(t)\mbox{Im}(\lambda_{j}(P))+\sum_{k=1}^{\infty} \frac{(-1)^k}{2k+1}[a(t)\mbox{Im}(\lambda_{j}(P))]^{2k+1},
\end{aligned}
\end{eqnarray}
so
\begin{eqnarray}\label{SAF_nn_9}
\prod_{t=0}^s \left(1-a(t)[1-\lambda_{j}(P)]\right)=\exp\Big(i\sum_{t=0}^{s} \varphi_t\Big) \prod_{t=0}^s r_t.
\end{eqnarray}
Assume $y_{j'}(0)\neq \frac{v_{j'}}{1-\lambda_{j}(P)}$.  Since
$\sum_{t=0}^{\infty}a(t)=\infty$, equations (\ref{SAF_nn_6}),
(\ref{SAF_nn_9}), and (\ref{SAF_nn_8}) imply
\begin{eqnarray}\label{SAF_nn_10}
\overline{\lim}_{s\to\infty}\overline{\lim}_{s_2\to\infty} |\E[y_{j'}(s_2)-y_{j'}(s)]|>0.
\end{eqnarray}
Next we consider the convergence of $x(s)$. Because $x(s)=H^{-1}
y(s)$, using Jensen's inequality we have
\begin{align}
\E\|x(s_2)-x(s)\|_2^2&=\E\|H^{-1}[y(s_2)-y(s)]\|_2^2 \nonumber\\
&\geq \sigma_n^2(H^{-1}) \E\|y(s_2)-y(s)\|_2^2 \nonumber\\
&\geq \sigma_n^2(H^{-1}) \E|y_{j'}(s_2)-y_{j'}(s)|^2 \nonumber\\
&\geq \sigma_n^2(H^{-1}) |\E[y_{j'}(s_2)-y_{j'}(s)]|^2\label{SAF_nn_11},
\end{align}
where $\sigma_n(H^{-1})=\inf_{\|x\|_2=1} \|H^{-1} x\|_2$ denotes the
least singular value of $H^{-1}$. Because $H^{-1}$ is invertible, we
have $\sigma_n(H^{-1})>0$. Hence, by (\ref{SAF_nn_10}) and
(\ref{SAF_nn_11}), we obtain
\begin{eqnarray*}\label{SAF_nn_12}
\overline{\lim}_{s\to\infty}\overline{\lim}_{s_2\to\infty} \E\|x(s_2)-x(s)\|_2^2>0.
\end{eqnarray*}
By the Cauchy criterion (see~\cite[page~58]{JAH:07}), $x(s)$ is not mean
square convergent.\\ (b) The geometric multiplicity of the eigenvalue
$1$ is less than its algebraic multiplicity.  By (a), we only need to
consider the case when any eigenvalue of $P$ with $1$ as real part has
zero imaginary part. Thus, the Jordan normal form $D$ contains a
Jordan block
\begin{eqnarray*}\label{SAF_nn_13}
J_j=
\begin{bmatrix}
    1 & 1 &  &     \\
     & 1 & \ddots &     \\
     &  & \ddots & 1   \\
     &  &  &   1
\end{bmatrix}_{m_j\times m_j}
\end{eqnarray*}
with $m_j\geq 2$.  Let $j'$ be the row index of $D$ corresponding to
the second line from the bottom of $J_{j}$.  It can be computed that
\begin{eqnarray*}\label{SAF_nn_14}
  [B(s)\cdots B(t)]_{j',j'+1}=\sum_{k=t}^s a(k).
\end{eqnarray*}
Since $\sum_{k=0}^{\infty}a(k)=\infty$, from (\ref{SAF_nn_2}), there
are some initial states such that $\lim_{s\to\infty}
|\E[y_{j'}(s)]|=\infty,$ which is followed by
$\lim_{s\to\infty}\E\|x(s)\|_2^2=\infty$.\\ (c) There is a left
eigenvector $\xi^T$ of $P$ corresponding to the eigenvalue
$1$ such that $\xi^T u\neq 0$. By (\ref{SAF_m1}) and (A1) we have
 \begin{align*}
\xi^T \E x(s+1)&=(1-a(s))\xi^T \E x(s)+a(s)[ \xi^T P \E x(s)+\xi^T u]\\
&=\xi^T \E x(s)+a(s)\xi^T u\\
&=\cdots=\xi^T x(0)+ \sum_{k=0}^s a(k) \xi^T u ,
\end{align*}
which implies $\lim_{s\to\infty} \E\|x(s)\|_2^2=\infty$ by $\sum_{k=0}^{\infty} a(k)=\infty$.\\
ii) It can be obtained by the similar method as i).

\section{Proof of Theorem \ref{SAF_r2}}\label{Proof_SAF_r2}
We prove our result by contradiction: Suppose that there exists a real
number sequence $\{a(s)\}_{s\geq 0}$ independent with $\{x(s)\}$ such
that
\begin{eqnarray}\label{SAF_r2_1}
\lim_{s\to\infty} \E\big\|x(s)-b\big\|_2^2=0.
\end{eqnarray}
We assert that $\lim_{s\to\infty} a(s)=0$. This assertion will be proved still by contradiction:
Assume that there exists a subsequence $\{a(s_k)\}_{k\geq 0}$ which does not converge to zero.
Let
$\widetilde{P}(s):=P(s)-P$ and $\widetilde{u}(s)=u(s)-u$ for any $s\geq 0$, then by (\ref{SAF_m1}), (A1) and (\ref{SAF_r2_01}) we have
\begin{align}
\E\Big[&\big\|x(s_k+1)-b\big\|_2^2\,|\,x(s_k)\Big]\nonumber\\
&=\E\Big[\big\|\xi+a(s_k)\big(\widetilde{P}(s_k) x(s_k)+\widetilde{u}(s_k) \big)\big\|_2^2\,|\,x(s_k)\Big]\nonumber\\
&=\|\xi\|_2^2 +a^2(s_k)\E\Big[\big\|\widetilde{P}(s_k) x(s_k)+\widetilde{u}(s_k)\big\|_2^2\,|\,x(s_k)\Big]\nonumber\\
&\geq a^2(s_k) c_3\label{SAF_r2_2},
\end{align}
where $$\xi:=(1-a(s_k))x(s_k)+a(s_k)(Px(s_k)+u)-b.$$ From
(\ref{SAF_r2_2}) we know that $\E\|x(s_k+1)-b\|_2^2$ will not converge
to $0$ as $k$ grows to infinity, which is in contradiction with
(\ref{SAF_r2_1}).


Since $x(0)\neq b$, to guarantee the convergence of $x(s)$, the gain
function $\{a(s)\}_{s\geq 0}$ must at least contain one non-zero element.  Also,
from (\ref{SAF_r2_2}), we can obtain that the number of the non-zero
elements in the sequence $\{a(s)\}_{s\geq 0}$ must be infinite.  Thus,
together with the assertion of $\lim_{s\to\infty} a(s)=0$,
there exists an integer $s^*> 0$ such that $a(s^*-1)\neq 0$,
$\{a(i)\}_{i=0}^{s^*-2}$ contains non-zero element, and
\begin{eqnarray}\label{SAF_r2_3_1}
2|a(s)(1-\mbox{Re}(\lambda_j(P))|<1, ~\forall s\geq s^*, 1\leq j\leq n.
\end{eqnarray}
Let $A(s):=(1-a(s))I_n+ a(s)P(s)$.  By (\ref{SAF_m1}) we have
\begin{eqnarray*}\label{SAF_r2_3_2}
x(s+1)=A(s) x(s)+a(s)u(s),~~s\geq s^*.
\end{eqnarray*}
By (A1), we obtain
\begin{align*}
\E[&x(s+1)\,|\,x(s^*)]-(I_n-P)^{-1} u\\
&=\big[I_n- a(s)(I_n-P)\big]  \E[x(s)\,|\,x(s^*)]\\
&\quad+a(s) u -(I_n-P)^{-1} u\nonumber\\
&=\big[I_n-a(s)(I_n-P)\big] \big(\E[x(s)\,|\,x(s^*)]-(I_n-P)^{-1} u\big)\nonumber\\
&=\cdots=\Big(\prod_{k=s^*}^s \E[A(k)]\Big) \big(x(s^*)-(I_n-P)^{-1} u\big),\nonumber
\end{align*}
which implies
\begin{align}\label{SAF_r2_4}
\E[x(s+1)|x(s^*)]&=H^{-1} \Big(\prod_{i=s^*}^s [I_n -a(i)(I_n-D)]\Big) \\
&\quad \times H \big(x(s^*)-(I_n-P)^{-1} u\big)+(I_n-P)^{-1} u\nonumber
\end{align}
from (\ref{Jordan}). Set
\begin{align}
z(s):= &\Big(\prod_{i=s^*}^s [I_n -a(i)(I_n-D)]\Big) H\nonumber\\
&\cdot\big(x(s^*)-(I_n-P)^{-1} u\big)+H(I_n-P)^{-1}u-H b.\label{SAF_r2_5}
\end{align}
Using Jensen's inequality and (\ref{SAF_r2_4}) we have
\begin{align}
\E\big[ \|x(s+1)-b\|_2^2 \,|\,x(s^*)\big]&\geq \big\| \E\big[(x(s+1)-b)\,|\,x(s^*)\big] \big\|_2^2 \nonumber\\
&= \big\| \E[(x(s+1)\,|\,x(s^*)]-b \big\|_2^2\label{SAF_r2_6}\\
&=\| H^{-1} z(s) \|_2^2 \geq \sigma_n^2(H^{-1}) \| z(s) \|_2^2\nonumber,
\end{align}
where   $\sigma_n(H^{-1})=\inf_{\|x\|_2=1} \|H^{-1} x\|_2$ denotes the  least singular value of $H^{-1}$. Because $H^{-1}$ is invertible,
$\sigma_n(H^{-1})>0$.
Define
\begin{equation}\label{SAF_r2_6_1}
w_j(s):=\prod_{i=s^*}^s (1-a(i)[1-\lambda_j(P)])
\end{equation}
and
\begin{equation}\label{SAF_r2_6_2}
M_{j}:=\prod_{i=s^*}^{\infty} [I_{m_j} -a(i)(I_{m_j}-J_{j})].
\end{equation}
We can compute that
\begin{align*}
|w_j(s)|^2&=\prod_{i=s^*}^s |1-a(i)[1-\lambda_j(P)]|^2\\
&=\prod_{i=s^*}^s \big\{1-2a(i)[1-\mbox{Re}(\lambda_j(P))]\\
&\quad+a^2(i)[1-2\mbox{Re}(\lambda_j(P))+|\lambda_j(P)|^2 \big\}.
\end{align*}
From this and (\ref{SAF_r2_3_1}) we have $w_j(s)\neq 0$ for any finite $s$. Also, if
$w_j(\infty)=0$, then $[1-\mbox{Re}(\lambda_j(P))]\sum_{i=s^*}^{\infty} a(i)=\infty$.
Hence, by (A4) or (A4'), there exists a Jordan block $J_{j_1}$ associated with the eigenvalue
$\lambda_{j_1'}(P)$ such that
$w_{j_1'}(\infty)\neq 0$ and (\ref{Assump_ag_1}) holds.
Because $M_{j_1}$ is an upper triangular matrix whose diagonal elements are all $w_{j_1'}(\infty)\neq 0$,
we can obtain the least singular value
\begin{equation*}\label{SAF_r2_9}
\sigma_{m_{j_1}}(M_{j_1})>0.
\end{equation*}
Also, by (\ref{SAF_r2_6}) and (\ref{Jordan_2}), we obtain
\begin{align}
\E\|x(\infty)-b\|_2^2&=\E\big[\E[ \|x(\infty)-b\|_2^2 \,|\,x(s^*)]\big]\nonumber\\
&\geq \sigma_n^2(H^{-1}) \E \| z(\infty) \|_2^2\nonumber\\
&\geq  \sigma_n^2(H^{-1}) \E \| z(\infty)-Ez (\infty) \|_2^2\nonumber\\
&\geq \sigma_n^2(H^{-1}) \E \big\|\widetilde{I}_{j_1}[ z(\infty)-Ez (\infty)] \big\|_2^2\nonumber\\
&= \sigma_n^2(H^{-1}) \E \Big\|\widetilde{I}_{j_1}\Big(\prod_{i=s^*}^{\infty} [I_n -a(i)(I_n-D)]\Big)\nonumber\\
&\quad\times H\big(x(s^*)- Ex(s^*)\big) \Big\|_2^2\nonumber\\
&= \sigma_n^2(H^{-1}) \E \big\|M_{j_1} \widetilde{I}_{j_1} H\big(x(s^*)- Ex(s^*)\big) \big\|_2^2\nonumber\\
&\geq \sigma_n^2(H^{-1}) \sigma_{m_{j_1}}^2(M_{j_1}) \nonumber\\
&\quad\times\E \big\|\widetilde{I}_{j_1} H\big(x(s^*)- Ex(s^*)\big) \big\|_2^2.\label{SAF_r2_10}
\end{align}
Using (\ref{SAF_m1}) and (\ref{Assump_ag_1}) we have
\begin{align}
 \E &\left\{ \big\| \widetilde{I}_{j_1} H (x(s^*)-\E x(s^*))\big\|_2^2  \,|\,x(s^*-1)\right\}\nonumber\\
 &=a^2(s^*-1) \E \Big\{\big\|\widetilde{I}_{j_1} H\widetilde{P}(s^*-1)x(s^*-1)\nonumber\\
 &\quad+\widetilde{I}_{j_1} H \widetilde{u}(s^*-1)  \big\|_2^2\,|\,x(s^*-1)\Big\} \nonumber\\
 &\geq a^2(s^*-1)\left(c_1 \|x(s^*-1)\|_2^2 + c_2\right).\label{SAF_r2_11}
\end{align}
Because $c_1$ and $c_2$ cannot be zero at the same time, we consider the case when $c_2>0$ first.
With the fact that $a(s^*-1)\neq 0$ and (\ref{SAF_r2_11}) we obtain
\begin{align*}
 \E &\big\| \widetilde{I}_{j_1} H (x(s^*)-\E x(s^*))\big\|_2^2\\
 &=\E  \left\{  \E \big\| \widetilde{I}_{j_1} H (x(s^*)-\E x(s^*))\big\|_2^2  \,|\,x(s^*-1)\right\}>0.
\end{align*}
Substituting this into (\ref{SAF_r2_10}) yields
$\E\|x(\infty)-b\|_2^2>0$, which is contradictory with
(\ref{SAF_r2_1}).

For the case when $c_1>0$, by (\ref{SAF_r2_10}) and (\ref{SAF_r2_11}),
we have
\begin{align}
\E&\|x(\infty)-b\|_2^2\nonumber\\
&\geq \sigma_n^2(H^{-1}) \sigma_{m_{j_1}}^2(M_{j_1}) a^2(s^*-1)\nonumber\\
&\quad\cdot \E \left(\big\|\widetilde{I}_{j_1} H\widetilde{P}(0)x(s^*-1) \big\|_2^2 \,|\,x(s^*-1)\right)\nonumber\\
&\geq \sigma_n^2(H^{-1}) \sigma_{m_{j_1}}^2(M_{j_1}) a^2(s^*-1) c_1 \E\|x(s^*-1)\|_2^2.\label{SAF_r2_13}
\end{align}
 Because $\{a(i)\}_{i=0}^{s^*-2}$ contains non-zero elements, we set $s'$ to be the biggest number such that $s'\leq s^*-2$
and $a(s')\neq 0$. By (\ref{SAF_r2_2}) we have
\begin{align*}
\E\big\|x(s^*-1)\big\|_2^2&=\E\big\|x(s'+1)\big\|_2^2 \nonumber\\
&\geq a^2(s')EE\big[\big\|\widetilde{P}(s') x(s')+\widetilde{u}(s')\big\|_2^2\,|\,x(s')\big]\nonumber\\
&\geq a^2(s') c_3>0.
\end{align*}
Substituting this into (\ref{SAF_r2_13}) we get $\E\|x(\infty)-b\|_2^2>0$, which is
 contradictory with (\ref{SAF_r2_1}).

\renewcommand{\thesection}{\Roman{section}}
\section{Proof of Theorem \ref{SAF_r3}}\label{Proof_SAF_r3}
\renewcommand{\thesection}{\arabic{section}}
Similar to the proof of Theorem \ref{SAF_r2} we prove our result by contradiction: Suppose that there exists a real number sequence $\{a(s)\}_{s\geq 0}$ independent with $\{x(s)\}$ such that (\ref{SAF_r2_1}) holds. Since $x(0)\neq b$, by (\ref{SAF_r2_1}) $\{a(s)\}_{s\geq 0}$ must contain non-zero elements.
We consider the following three cases respectively to deduce the contradiction:\\
\textbf{Case I}: The condition i) is satisfied. Similar to the proof of Theorem \ref{SAF_r2}, we
first prove $\lim_{s\to\infty} a(s)=0$ by contradiction: Suppose there exists a subsequence $\{a(s_k)\}$ that does not converge to zero.
For the case when $b\neq {\mathbf{0}}_{n\times 1}$, by (\ref{SAF_r2_1}), there exists a time $s_1\geq 0$ such that
\begin{equation}\label{SAF_r3_1}
\E\big\|x(s)-b\big\|_2^2 \leq \frac{1}{4}\|b \|_2,~~\forall s>s_1.
\end{equation}
Because for any $x(s_k)$,
\begin{align*}
 \|b\|_2^2 & \leq \big(\|x(s_k)\|_2+ \big\|b-x(s_k)\big\|_2\big)^2\\
& \leq 2\big(\|x(s_k)\|_2^2+ \big\|b-x(s_k)\big\|_2^2\big),
\end{align*}
(\ref{SAF_r3_1}) is followed by
\begin{equation}\label{SAF_r3_3}
\E\|x(s_k)\|_2^2 \geq \frac{1}{2}\|b\|_2- \E\big\|x(s_k)-b\big\|_2^2 \geq \frac{1}{4}\|b\|_2
\end{equation}
for large $k$. By (\ref{SAF_r2_2}), (\ref{SAF_r3_01}) and (\ref{SAF_r3_3}) we obtain
\begin{align}
\E\big\|x(s_k+1)-b\big\|_2^2 &\geq a^2(s_k) \E\| \widetilde{P}(s_k) x(s_k)\|_2^2\nonumber\\
&= a^2(s_k)  \E\|H^{-1} H \widetilde{P}(s_k) x(s_k)\|_2^2\nonumber\\
&\geq a^2(s_k) \sigma_n^2(H^{-1}) \E\| H \widetilde{P}(s_k) x(s_k)\|_2^2\nonumber\\
&\geq a^2(s_k) \sigma_n^2(H^{-1}) c_1 \E\|x(s_k)\|_2^2\nonumber\\
&\geq \frac{1}{4}a^2(s_k) \sigma_n^2(H) c_1\|b\|_2,\label{SAF_r3_4}
\end{align}
which is contradictory with (\ref{SAF_r2_1}).

 For the case when $b={\mathbf{0}}_{n\times 1}$, by (\ref{SAF_r2_2}) and (\ref{SAF_r3_4}), we have
 \begin{small}
\begin{align}
\E\big\|x(s_k+1)\big\|_2^2&\geq \E\|(1-a(s_k))x(s_k)+a(s_k)(Px(s_k)+u)\|_2^2 \nonumber\\
&\quad+a^2(s_k) \sigma_n^2(H^{-1}) c_1 \E\|x(s_k)\|_2^2\label{SAF_r3_5}.
\end{align}
\end{small}
If $\|(1-a(s_k))I_n+a(s_k)P\|_2 \E\|x(s_k)\|_2 > \frac{1}{2}\|a(s_k)u\|_2$, by (\ref{SAF_r3_5}) and Jensen's inequality we have
\begin{small}
\begin{align}\label{SAF_r3_6}
\E\big\|x(s_k+1)\big\|_2^2&\geq a^2(s_k) \sigma_n^2(H^{-1}) c_1 (\E\|x(s_k)\|_2)^2\nonumber\\
&\geq \frac{a^4(s_k) \sigma_n^2(H^{-1}) c_1 \|u\|^2}{4 \|(1-a(s_k))I_n+a(s_k)P\|_2^2}\nonumber\\
& \nrightarrow 0 ~\mbox{if}~ a(s_k)\nrightarrow 0.
\end{align}
\end{small}
Otherwise,
\begin{align*}
\E&\|(1-a(s_k))x(s_k)+a(s_k)(Px(s_k)+u)\|_2\nonumber\\
&\geq \|a(s_k)u\|_2-\E\|(1-a(s_k))x(s_k)+a(s_k)Px(s_k)\|_2 \nonumber\\
&\geq  \|a(s_k)u\|_2-\E\|(1-a(s_k))I_n+a(s_k)P\|_2 \|x(s_k)\|_2\nonumber\\
&\geq \frac{1}{2}\|a(s_k)u\|_2,
\end{align*}
Hence, using (\ref{SAF_r3_5}) and Jensen's inequality again, we obtain
\begin{small}
\begin{align}\label{SAF_r3_8}
\E\big\|x(s_k+1)\big\|_2^2&\geq \big(\E\|(1-a(s_k))x(s_k)+a(s_k)(Px(s_k)+u)\|_2\big)^2 \nonumber\\
&\geq \|a(s_k)u\|_2^2/4.
\end{align}
\end{small}
Combining (\ref{SAF_r3_6}) and (\ref{SAF_r3_8}) yields
$\E\big\|x(s_k+1)\big\|_2^2$. This quantity does not converge to zero,
which is in contradiction with (\ref{SAF_r2_1}).  By summarizing the
arguments above we prove the assertion of
$\lim_{s\to\infty}a(s)=0$.

Because $\lim_{s\to\infty}a(s)=0$ and because $\{a(s)\}_{s\geq 0}$
contains non-zero elements, there exists an integer $s^*>0$ such that
$a(s^*-1)\neq 0$ and (\ref{SAF_r2_3_1}) holds.  Define $w_j(s)$ and
$M_{j}$ by (\ref{SAF_r2_6_1}) and (\ref{SAF_r2_6_2}) respectively.
With the arguments similar to the proof of Theorem \ref{SAF_r2}, we can
find a Jordan block $J_{j_1}$ associated with the eigenvalue
$\lambda_{j_1'}(P)$ such that $w_{j_1'}(\infty)\neq 0$ and
(\ref{SAF_r3_01}) holds. Similar to (\ref{SAF_r2_13}) we obtain
\begin{small}
\begin{equation}\label{SAF_r3_9}
\E\|x(\infty)-b\|_2^2\geq \sigma_n^2(H) \sigma_{m_{j_1}}^2(M_{j_1})  a^2(s^*-1) c_1 \E\|x(s^*-1)\|_2^2.
\end{equation}
\end{small}
By (\ref{SAF_r3_5}) we have that if $\E\|x(s)\|_2^2>0$, then
$\E\|x(s+1)\|_2^2>0$ for any $a(s)\in\real$.  Then with the condition
$x(0)\neq {\mathbf{0}}_{n\times 1}$, we have $\E\|x(s^*-1)\|_2^2>0$.  Using this and
(\ref{SAF_r3_9}) we get $\E\|x(\infty)-b\|_2^2>0$, which is
contradictory with (\ref{SAF_r2_1}).

\textbf{Case II}: The condition ii) is satisfied. Since
$\{a(s)\}_{s\geq 0}$ contains non-zero elements, we define $s_1$ to be
the first $s$ such that $a(s)\neq 0$. Then $x(s_1+1)=a(s_1)u\neq b$
almost surely. Let $s_1+1$ be the initial time and by the same
arguments as in Case I we obtain $\E\|x(\infty)-b\|_2^2>0$.

\textbf{Case III}: The condition iii) is satisfied. If $x(0)={\mathbf{0}}_{n\times 1}$,
we obtain $\E\|x(s)\|_2^2=0$ for any $s\geq 0$, which is contradictory
with (\ref{SAF_r2_1}). Thus, we just need to consider the case when
$x(0)\neq {\mathbf{0}}_{n\times 1}$. Since $\{a(s)\}_{s\geq 0}$ contains non-zero
elements, we define $s_1$ to be the first $s$ such that $a(s)\neq 0$.

Set $x^*:=x_1+1$. Define $w_j(s)$ and $M_{j}$ by (\ref{SAF_r2_6_1})
and (\ref{SAF_r2_6_2}) respectively.  If $\lambda_j(P)$ is not a real
number, then $w_{j}(s)$ cannot be equal to $0$ for any finite $s$. By the
similar arguments as in the proof of Theorem \ref{SAF_r2}, there exists a
Jordan block $J_{j_1}$ associated with the eigenvalue
$\lambda_{j_1'}(P)$ such that $w_{j_1'}(\infty)\neq 0$ and
(\ref{SAF_r3_01}) holds. By (\ref{SAF_r3_9}) we have
\begin{small}
\begin{align*}
\E\|x(\infty)-b\|_2^2
&\geq \sigma_n^2(H^{-1}) \sigma_{m_{j_1}}^2(M_{j_1}) a^2(s^*-1) c \E\|x(s^*-1)\|_2^2\\
&=\sigma_n^2(H^{-1}) \sigma_{m_{j_1}}^2(M_{j_1}) a^2(s_1) c \|x(0)\|_2^2>0,
\end{align*}
\end{small}
which is contradictory with (\ref{SAF_r2_1}).

\ifCLASSOPTIONcaptionsoff
  \newpage
\fi

\end{document}